\documentclass[10pt, a4paper,draft]{amsart}

\usepackage{amsmath} \usepackage{amssymb} \usepackage{amsfonts}

\DeclareMathOperator{\spt}{supp}

\numberwithin{equation}{section}
\newtheorem{theorem}{Theorem}[section]
\newtheorem{corollary}[theorem]{Corollary}
\newtheorem{proposition}[theorem]{Proposition}
\newtheorem{lemma}[theorem]{Lemma}

\theoremstyle{remark}

\theoremstyle{definition} \newtheorem{definition}{Definition}[section]

\theoremstyle{remark} 
\theoremstyle{remark} \newtheorem{remark}{Remark}[section]

\begin{document}

\newcommand{\diff}{{\rm\, d}} \newcommand{\R}{{\mathbb R}}
\newcommand{\C}{{\mathbb C}} \newcommand{\D}{\mathcal D}
\newcommand{\xla}{\langle\lambda\rangle}
\newcommand{\xxi}{\langle\xi\rangle} \newcommand{\xx}{\langle  x\rangle}
\newcommand{\yy}{\langle y\rangle}
\newcommand{\DD}{\langle D\rangle}
\newcommand{\dint}{\int\!\!\int}
\newcommand{\triple}[1]{{|\!|\!|#1|\!|\!|}} \newcommand{\Sph}{{\mathbb
    S}} \newcommand{\Ls}{{L^{2}_{s}}}

\title[Wave and Dirac equations with magnetic
potential]%
{Decay estimates for the wave and Dirac equations with a magnetic
 potential}

\begin{abstract}
   We study the electromagnetic wave equation
   and the perturbed massless Dirac equation
   on $\mathbb{R}_{t}\times\mathbb{R}^{3}$:
   \begin{equation*}
    u_{tt}-(\nabla+iA(x))^{2}u+B(x)u=0,\qquad
    iu_t-\D u+V(x)u=0
  \end{equation*}
  where the potentials $A(x),B(x),V(x)$ are assumed to be
  small but may be rough.
  For both equations, we prove
  the expected time decay rate of the solution
  \begin{equation*}
    |u(t,x)|\leq\frac{1}{t}\|f\|_X
  \end{equation*}
  where the norm $\|f\|_{X}$
  can be expressed as the weighted $L^{2}$
  norm of a few derivatives
  of the data $f$.
  \end{abstract}

\date{}    

\author{Piero D'Ancona} \address{Piero D'Ancona: Unversit\`a di Roma
  ``La Sapienza'', Dipartimento di Matematica, Piazzale A.~Moro 2,
  I-00185 Roma, Italy} \email{dancona@mat.uniroma1.it}

\author{Luca Fanelli} \address{Luca Fanelli: Unversit\`a di Roma ``La
  Sapienza'', Dipartimento di Matematica, Piazzale A.~Moro 2, I-00185
  Roma, Italy} \email{fanelli@mat.uniroma1.it}

 \maketitle

\section{Introduction}\label{sec.introd}

Dispersive properties of evolution equations play a
crucial role in the study of nonlinear problems, and for this
reason they have attracted a great deal of attention in
recent years. In particular, for the Schr\"odinger and the
wave equation a well established theory exists, see
\cite{GinibreVelo95-generstric}
and
\cite{KeelTao98-endpoinstric}.
On the other hand, in the variable coefficient
case the theory is very far from complete.
The simplest situtation is a perturbation with a term of order
zero; this is already very interesting from the physical point of view
(electrostatic potential). Several results are available for the
equations
\begin{equation*}
    i\partial_{t}u-\Delta u+V(x)u=0,\qquad\qquad
    \square u +V(x)u=0.
\end{equation*}
We cite among the others
\cite{BurqPlanchonStalkerTahvildar03},
\cite{Goldberg04},
\cite{GoldbergSchlag04},
\cite{JourneSofferSogge91-decayschroed},
\cite{RodnianskiSchlag04-disp}
and the recent survey \cite{Schlag05}
for Schr\"o\-ding\-er; and
\cite{Beals94-optiml},
\cite{BealsStrauss93-l},
\cite{Cuccagna00},
\cite{pda-vf},
\cite{GeorgievVisciglia03}
for the wave equation. We must also
mention the wave operator approach of Yajima (see
\cite{ArbYaj00},
\cite{Yajima95-waveopN},
\cite{Yajima95-waveopNeven},
\cite{Yajima99-waveop2D})
which permits to deal with the above equations
in a unified way, although under nonoptimal assumptions on the
potential in dimensions 1 and 3.

The next step in generality is a first order perturbation; from the
physical point of view this corresponds to a magnetic potential. In
this case only a handful of results are available: Strichartz
estimates for the 3D wave equation \cite{CuccagnaSchirmer01},
provided the coefficients are small and in the Schwartz class; and
smoothing estimates for the 3D Schr\"odinger and wave operators
\cite{Tarulli04-smoo}. The most general case of variable
coefficients has been studied in
\cite{HassellTaoWunsch04},
\cite{RobbianoZuily05} and
\cite{StaffilaniTataru02-stricschroed},
where local Strichartz estimates have been proved, in various
degrees of complexity; see also \cite{Burq04}.

In the present paper, our main focus will be on the three dimensional wave
equation with an electromagnetic potential
\begin{equation}\label{eq.wmag}
    u_{tt}-(\nabla +iA(x))^{2}u+B(x)u=0,\qquad
    u:\mathbb{R}\times\mathbb{R}^{3}\to\mathbb{C},
\end{equation}
 and the closely related massless Dirac system with a potential:
 \begin{equation}\label{eq.Dpot}
    iu_{t}-\D u+V(x)u=0,\qquad
    u:\mathbb{R}\times\mathbb{R}^{3}\to\mathbb{C}^{4}.
\end{equation}
Here $A:\mathbb{R}^{3}\to\mathbb{R}^{3}$,
$B:\mathbb{R}^{3}\to\mathbb{R}$, $V(x)=V^{*}(x)$ is
a 4$\times$4 complex matrix on $\mathbb{R}^{3}$,
and the symbol $\D$ denotes the constant coefficient,
elliptic, $L^{2}$ selfadjoint operator
\begin{equation*}
  \D=\frac 1 i \sum_{j=1}^{3}\alpha_{k}{\partial_{k}},
\end{equation*}
where the \emph{Dirac matrices} $\alpha_{1}$, $\alpha_{2}$,
$\alpha_{3}$ have the following structure:
\begin{equation}\label{eq.diracm}
  \alpha_{1}=
  \begin{pmatrix}
    0 & 0 & 0 & 1 \\
    0 & 0 & 1 & 0 \\
    0 & 1 & 0 & 0 \\
    1 & 0 & 0 & 0
  \end{pmatrix},\quad \alpha_{2}=
  \begin{pmatrix}
    0 & 0 & 0 & -i \\
    0 & 0 & i & 0 \\
    0 & -i & 0 & 0 \\
    i & 0 & 0 & 0
  \end{pmatrix},\quad \alpha_{3}=
  \begin{pmatrix}
    0 & 0 & 1 & 0 \\
    0 & 0 & 0 & -1 \\
    1 & 0 & 0 & 0 \\
    0 & -1 & 0 & 0
  \end{pmatrix}.
\end{equation}
We neglect the physical constants (i.e., we set $c=\hbar =1$),
and we consider the zero mass case exclusively; the case of
a positive mass, whose
second order counterpart is the Klein-Gordon equation,
has an additional term $\alpha_{4}u$ with
\begin{equation}\label{eq.alpha4}
  \alpha_{4}=
  \begin{pmatrix}
    1 & 0 & 0 & 0 \\
    0 & 1 & 0 & 0 \\
    0 & 0 & -1 & 0 \\
    0 & 0 & 0 & -1
  \end{pmatrix}.
\end{equation}
The relation between massless Dirac and wave equation
is readily explained: indeed, the Dirac matrices satisfy the
commutation rules
\begin{equation*}
\alpha_{\ell}\alpha_{k}+\alpha_{k}\alpha_{\ell}
=2\delta_{kl}I_{4}
\end{equation*}
which imply immediately
\begin{equation*}
     \D^{2}=-\Delta I_{4},
\end{equation*}
where $I_{4}$ is the 4$\times$4 identity matrix. Thus
we have the fundamental relation
\begin{equation*}
  (i\partial_{t}-\D)(i\partial_{t}+\D)=
  (\Delta-\partial^{2}_{tt}) I_{4},
\end{equation*}
which can be intepreted as follows: squaring the
Dirac system produces a diagonal system of wave equations
(or, conversely: taking the square root of a wave equation
produces a Dirac system. According to the folklore, this
was the route that lead Dirac to his equation).
When a potential is present in the Dirac system, the above reduction
produces an electromagnetic wave equation in a natural way.
A discussion of this can be found e.g. in
\cite{LandauLifshitz}
(Volume 4, Chapter 4); see also section \ref{sec.th2} below.

%
%
%

Our goal here is to establish the decay rate of the spatial
$L^{\infty}$ norm
of the solution, with minimal assumptions on the potentials.
The expected decay rate  is $t^{-1}$, both for the
wave equation and the Dirac system. Indeed, known results
for hyperbolic systems (for constant coefficients see e.g.
\cite{Liess91-crystal},
\cite{LucenteZiliotti00-maxw},
and for $C^{\infty}_{0}$ perturbations thereof see
\cite{KajitaniSatoh04})
suggest a
$t^{-\frac{n-1}2}$ decay rate
in $n$ space dimensions.

Before stating our first result we introduce some
basic notations.
Under the assumptions of Theorem \ref{th.1} below,
the perturbed laplacian
\begin{equation}\label{eq.magnlap}
    H:=-(\nabla+iA(x))^2+B(x),
\end{equation}
where $A(x)=(A_1(x),A_2(x),A_3(x)):
\mathbb{R}^{3}\to\mathbb R^3$ and
$B(x):\mathbb{R}^{3}\to\mathbb R$, is a selfadjoint
unbounded operator on $\mathbb{R}^{3}$; the explicit
standard construction is recalled in Section \ref{sec.self}.
Spectral calculus allows us to define the operators
$\psi(H)$ for any well behaved function $\psi(s)$.

In particular, consider a \emph{(non-homogeneous)
Paley-Littlewood partition of unity} on $\mathbb{R}^{3}$,
defined as follows:
fix a radial nonnegative function $\psi(r)\in C^{\infty}_{0}$
with $\psi(r)=1$ for $r<1$, $\psi(r)=0$ for $r>2$, define
$\phi_{j}(r)=\psi(2^{-j+2}r)-\psi(2^{-j+1}r)$ for all $j\geq1$,
and $\phi_{0}=\psi$. Then $1=\sum_{j\geq0}\phi_{j}$ is
the required partition of unity on $\mathbb{R}^{3}$. The
operators $\phi_{j}(\sqrt H)$ will be used in the following
to define suitable norms associated to the operator $H$.
We shall also use the notations
\begin{equation*}
    \xx=(1+|x|^{2})^{1/2},\qquad
    \DD^{s} f=(1-\Delta)^{s/2}f
         \equiv \mathcal{F}^{-1}(\xxi^{s}\widehat f)
\end{equation*}

Our first result concerns the Cauchy problem for the
wave equation perturbed with a small rough electromagnetic
potential
\begin{align}\label{eq.onde}
  &u_{tt}(t,x)-(\nabla+iA(x))^2u+B(x)u=0,\qquad
     (t,x)\in\mathbb{R}\times\mathbb{R}^{3},   \\
  &u(0,x)=0,\quad u_t(0,x)=g(x).\label{eq.onded}
\end{align}
We can prove:

\begin{theorem}\label{th.1}
  Assume the potentials
  $A(x)\in\mathbb{R}^{3}$, $B(x)\in\mathbb{R}$ satisfy
  \begin{equation}\label{eq.hpWfin}
    |A_{j}|\leq \frac{C_0}{|x|\xx(|\log|x||+1)^\beta},\qquad
    \sum_{j=1}^3|\partial_jA_j|+
    |B|\leq\frac{C_0}{|x|^2(|\log|x||+1)^\beta},
  \end{equation}
  for some constant $C_{0}>0$ sufficiently small and some
  $\beta>1$. Then any
  solution of the Cauchy problem \eqref{eq.onde},
  \eqref{eq.onded} satisfies the decay estimate
  \begin{equation}\label{eq.decwave}
    |u(t,x)|\leq
    \frac{C}{t}\sum_{j\geq0}2^{2j}
        \|\xx w_\beta^{1/2}
        \varphi_j(\sqrt{H})g\|_{L^2},
  \end{equation}
  where $w_\beta(x):=|x|(|\log|x||+1)^\beta$.
  If in addition we assume that, for some $\epsilon>0$,
  \begin{equation}\label{eq.ipfinregW}
      \DD^{1+\epsilon}A_j \in L^\infty,\qquad
   \DD^{\epsilon}B\in {L^\infty}
  \end{equation}
  then $u$ satisfies  for any $\delta>0$
   the estimate
  \begin{equation}\label{eq.decwaveH2}
    |u(t,x)|\leq
    \frac{C}{t}\|\xx^{3/2+\delta}g\|_{H^{2+\epsilon}}.
  \end{equation}
\end{theorem}

\begin{remark}\label{rem.norm}
The norm appearing in \eqref{eq.decwave} can be regarded as a
distorted analogue of a standard Besov norm, generated by the
operator $H$. Similar norms already appeared in
\cite{CuccagnaSchirmer01} for magnetic potentials with coefficients
in the Schwartz class; in that case, however, it was possible to
prove the equivalence with standard Besov norms (see also
\cite{pda-vf}, \cite{GeorgievVisciglia03} for the analogous norms
generated by $-\Delta+V(x)$, which are also equivalent to the
nondistorted norms). Under the slightly stronger assumptions
\eqref{eq.ipfinregW} on the coefficients, it is possible to prove an
estimate like \eqref{eq.decwaveH2} expressed in terms of standard
weighted Sobolev norms.

Moreover, we remark that in
our estimates we lose 2 derivatives; it is natural to conjecture that
this is not optimal, and it should be possible to lose only one derivative
as in the case of the free wave equation.
\end{remark}

\begin{remark}\label{rem.psu}
As an essential step in the proof of Theorem \ref{th.1}, we need to
establish the \emph{limiting absorption principle} (LAP) for the
operator $H$. This is obtained in Section \ref{sec.lap} through
several steps: starting from the ``weak''  LAP of \cite{BRV} for the
free resolvent, we first prove a strong version of the LAP for the
free operator in the weighted spaces
\begin{equation*}
    L^{2}(w_{\beta}(x)dx),\qquad
    w_\beta(x):=|x|(|\log|x||+1)^\beta
\end{equation*}
and then we get the LAP  for the perturbed operator. For the precise
statements see Proposition \ref{prop.lapW}. See also
\cite{Tarulli04-smoo} for related results.
\end{remark}

\begin{remark}\label{rem.fg}
When the initial data are of the form
\begin{equation*}
    u(0,x)=f,\qquad u_{t}(0,x)=0,
\end{equation*}
Theorem \ref{th.1} implies, by standard arguments,
the estimate
\begin{equation}\label{eq.decwavef}
    |u(t,x)|\leq
    \frac{C}{t}\sum_{j\geq0}2^{3j}
        \|\xx w_\beta^{1/2}\varphi_j(\sqrt H)f\|_{L^2}
\end{equation}
with an additional loss of one derivatives as expected.
If in addition we assume that for some $\epsilon>0$
\begin{equation}\label{eq.ipfinregWf}
      \DD^{2+\epsilon}A_j \in L^\infty,\qquad
   \DD^{1+\epsilon}B\in {L^\infty}
\end{equation}
then also the simpler estimate
\begin{equation}\label{eq.decwaveH2f}
    |u(t,x)|\leq
    \frac{C}{t}\|\xx^{3/2+\delta}f\|_{H^{3+\epsilon}}.
\end{equation}
holds for all $\delta>0$.
\end{remark}

Our second result concerns the perturbed Dirac system
\begin{align}\label{eq.Dirac}
  &iu_{t}-\D u+V(x)u=0,\qquad
  (t,x)\in\mathbb{R}\times\mathbb{R}^{3},\\
  &u(0,x)=f(x).\label{eq.Diracd}
\end{align}
By explointing the above mentioned
relation between the magnetic wave equation
and the Dirac system, we can prove the following Theorem
as a direct consequence of Theorem \ref{th.1}:

\begin{theorem}\label{th.2}
  Assume the 4$\times$4 complex valued matrix
  $V(x)=V^{*}(x)$ satisfies
    \begin{equation}\label{ipfinVreg}
        |V(x)|\leq\frac{C_0}{|x|\xx(|\log|x||+1)^\beta},\qquad
        |DV(x)|\leq\frac{C_0}{|x|^2(|\log|x||+1)^\beta},
    \end{equation}
    for some $C_0>0$ small enough
    and some $\beta>1$.
    Then the solution of the Cauchy problem \eqref{eq.Dirac},
    \eqref{eq.Diracd} satisfies the decay estimate
    \begin{equation}\label{estlastdir}
        |u(t,x)|\leq\frac Ct\sum_{j\geq0}2^{3j}\|\xx w_\beta^{1/2}
        \varphi_j(\D+V)f\|_{L^2},
    \end{equation}
    where $w_\beta(x)=|x| (|\log|x||+1)^\beta$.
    If in addition we assume that, for some $\epsilon>0$,
    \begin{equation}\label{ipfinVregH3}
        \DD^{2+\epsilon}V\in L^\infty,
    \end{equation}
  then $u$ satisfies  for any $\delta>0$
   the estimate
    \begin{equation}\label{estlastdirH3}
        |u(t,x)|\leq\frac Ct\|\xx^{3/2+\delta}f\|_{H^{3+\epsilon}}.
    \end{equation}
\end{theorem}

Since Theorem \ref{th.2} is proved essentially by
``squaring'' the perturbed Dirac operator, a condition
on the derivative $DV$ is essential in order to apply
Theorem \ref{th.1} to the resulting wave equation.
On the other hand, we can study the Cauchy problem
\eqref{eq.Dirac}, \eqref{eq.Diracd} by a direct
application of the spectral calculus for the selfadjoint
operator $\D+V(x)$; this alternative approach allows
us to consider much rougher potentials $V(x)$
(see \eqref{eq.ipfinV}). The price to pay is an
additional loss of one derivative, so that the total
loss is 4 derivatives in our last result:

\begin{theorem}\label{th.3}
  Assume the 4$\times$4 complex valued matrix
  $V(x)=V^{*}(x)$ satisfies
  \begin{equation}\label{eq.ipfinV}
    |V(x)|\leq
    \frac{C_0}{|x|^{1/2}\xx^{3/2}(|\log|x||+1)^{\beta/2}},
  \end{equation}
  for some $C_0>0$ small enough and some
  $\beta>1$. Then the solution of the Cauchy problem
  \eqref{eq.Dirac},  \eqref{eq.Diracd}
  satisfies for any $\epsilon>0$ the decay estimate
  \begin{equation}\label{eq.decayVfin}
    |u(t,x)|\leq\frac Ct\sum_{j\geq0}2^{4j}
    \|\xx^{3/2+\epsilon}\varphi_j(\D+V)f\|_{L^2}.
  \end{equation}
\end{theorem}

\begin{remark}\label{rem.lapdir}
As a byproduct of our method of proof, we obtain the \emph{limiting
absorption principle} for the perturbed Dirac operator under
assumption \eqref{eq.ipfinV} (see Section \ref{ssec.dirac}). The LAP
had been proved earlier for the free Dirac equation by Yamada
\cite{yamada}, and for the Dirac equation with potential (and with
mass) in \cite{PSU} under quite stronger assumptions.
\end{remark}

\section{The self-adjointness of the perturbed operators}
\label{sec.self}

In this section we check the selfadjointness of
the perturbed operators $\Delta_W$ and $\mathcal
D_V$ under quite general assumptions on the potentials
$A,B,V$, which in particular are implied by the
assumptions of Theorems \ref{th.1}, \ref{th.2}
and \ref{th.3}. Most of the material here is standard;
however we decided to include a sketch of the proof
for the sake of completeness. Moreover, the use
of Lorentz spaces techniques (see the Appendix
for a short review) makes the proofs quite straightforward.


It will be useful sometimes to express the magnetic laplacian both in
the covariant form
\begin{equation}\label{eq.cov}
    H=-(\nabla+iA(x))^{2}+B(x)
\end{equation}
and in the expanded form
\begin{equation}\label{eq.Wdoppio}
H=-\Delta+W(x,D),\qquad
W(x,D)=\sum_{j=1}^3a_j(x)\partial_j+b(x)
\end{equation}
where
\begin{equation}\label{eq.coeffi}
    a_j(x)=-2iA_j(x), \qquad b(x)=-i\sum_{j=1}^3\partial_j A_j(x)
    +|A(x)|^2+B(x),\qquad A_{j},B\in\mathbb{R}.
\end{equation}
Then we have the following:

\begin{proposition}\label{prop.selfadjW}
    Consider the operator on $C^{\infty}_{0}(\mathbb{R}^{n})$
    \begin{equation}\label{eq.magnlapn}
        H=-(\nabla+iA(x))^{2}+B(x),
    \end{equation}
    where $A(x):\mathbb{R}^{n}\to\mathbb{R}^{n}$ and
    $B(x):\mathbb{R}^{n}\to\mathbb{R}$ are measurable functions.
    Assume that the Lorentz (weak Lebesgue) norms
    of the coefficients
    \begin{equation}\label{eq.selfadj}
        \|A\|_{L^{n,\infty}}\leq C_0,\qquad
        \|B\|_{L^{n/2,\infty}}\leq C_0
    \end{equation}
    are bounded by
    some constant $C_0>0$ small enough. Then
    $H$ has a (unique) self-adjoint extension
    to $H^2(\mathbb R^n)$.
\end{proposition}

\begin{proof}
    Our proof is based on the standard results on
    quadratic forms, see e.g. the standard reference
    \cite{reed-simonI}. First of
    all we notice that by \eqref{eq.selfadj} we have immediately
    \begin{equation*}
    |A(x)|^{2}\in L^{n/2,\infty}
    \end{equation*}
    with a small norm.
    Now, the quadratic form $q(\phi,\psi)$ given by
    $$q(\varphi,\psi)
    =((\nabla+iA)\varphi,(\nabla+iA)\psi)_{L^{2}}
    +(B\varphi,\psi)_{L^{2}}$$
     is well defined on the form domain  $H^{1}$ under assumptions
     \eqref{eq.selfadj}. Indeed,  consider the identity
     \begin{equation}\label{eq.identq}
    q(\psi,\psi)=\|\nabla\psi\|_{L^{2}}^{2}+
       ((|A|^{2}+B)\psi,\psi)_{L^{2}}
       +2\Im(A\nabla\psi,\psi)_{L^{2}};
    \end{equation}
    using the embedding $H^1\subset L^{2n/(n-2),2}$,
    the H\"older inequality in Lorentz spaces
    (see the Appendix at the end of the paper for a quick synopsis of the relevant
    results), and recalling assumption \eqref{eq.selfadj}, we have easily
    \begin{align*}
    |q(\psi,\psi)|\leq &\|\nabla\psi\|_{L^{2}}^{2}+
           C\||A|^{2}+B\|_{L^{n/2,\infty}}
                      \|\psi\overline\psi\|_{L^{\frac n{n-2},1}}+
           C\|A\|_{L^{n,\infty}}
                     \|\nabla\psi\cdot
                          \overline\psi\|_{L^{\frac n{n-1},1}}\\
       \leq&\|\nabla\psi\|_{L^{2}}^{2}+
           CC_{0}
                      \|\psi\|_{L^{\frac {2n}{n-2},2}}^{2}+
           CC_{0}
                     \|\nabla\psi\|_{L^{2,2}}
                          \|\psi\|_{L^{\frac {2n}{n-2},2}}
           \leq C\|\nabla\psi\|_{L^{2}}^{2}.
   \end{align*}
    It is clear that the form is symmetric, since
    $A$ and $B$ are real valued.
    Now, recalling Theorem VIII.15  in \cite{reed-simonI}, in order
    to prove that $q$ is the form associated to a
    (uniquely defined) self-adjoint operator,
    it will be sufficient to show that it is \emph{closed}, i.e., its
    domain $H^{1}(\mathbb{R}^{n})$ is complete under the
    norm
    \begin{equation}\label{eq.normq}
    |\!|\!|\psi|\!|\!|^{2}=q(\psi,\psi)
       +C\|\psi\|_{L^{2}}^{2}
    \end{equation}
    for some $C>0$,
    and that it is \emph{semibounded}, i.e.,
    \begin{equation}\label{eq.semib}
    q(\psi,\psi)\geq-C\|\psi\|_{L^{2}}^{2}
   \end{equation}
   for some $C>0$.
   Both properties follow from the identity
   \eqref{eq.identq}; indeed, by
   estimating as above we obtain easily
   \begin{equation*}
    q(\psi,\psi)\geq \|\nabla\psi\|^{2}_{L^{2}}-
       CC_{0}\|\nabla\psi\|_{L^{2}}^{2}.
    \end{equation*}
    In particular this implies that the norm \eqref{eq.normq}
    is \emph{equivalent} to the $H^{1}(\mathbb{R}^{n})$
    norm, provided $C_{0}$ is small enough, so that the form is
    closed; and this implies also that \eqref{eq.semib} is satisfied
    with $C=0$.
\end{proof}

For the perturbed Dirac operator we have a similar result:

\begin{proposition}\label{prop.selfadjV}
    Let $V(x)=V^{*}(x)$ be a 4$\times$4 complex valued
    matrix on $\mathbb R^3$.
    Assume that
    \begin{equation}\label{eq.selfadjV}
        \|V\|_{L^{3,\infty}}\leq C_0,
    \end{equation}
    for some $C_0>0$ sufficiently small. Then the perturbed
    Dirac operator $\D_V=\D+V$ is self-adjoint on $H^1(\mathbb
    R^3,\mathbb C^4)$.
\end{proposition}

\begin{proof}
    The proof is analogous to the proof of Theorem
    \ref{prop.selfadjW}. We
    define the quadratic form $q:H^{1/2}\times H^{1/2}\to \mathbb C$
    associated to the operator $\D_V$ as
    $$q(\varphi,\psi):=(\D\varphi,\psi)+(V\varphi,\psi).$$
    First we prove that the domain of $q$ is $H^{1/2}$.
    With the same arguments of the previous theorem we
    estimate
    \begin{eqnarray*}
        |q(\varphi,\varphi)| & \leq & \|\varphi\|_{H^{1/2}}^2
        +C\|V\|_{L^{3,\infty}}\|\varphi^2\|_{L^{n/(n-1),1}}
        \\ \  & \leq & \|\varphi\|_{H^{1/2}}^2
        +C\|V\|_{L^{3,\infty}}\|\varphi\|_{L^{2n/(n-1),2}}^2
        \\ \  & \leq & \left(1+C\|V\|_{L^{3,\infty}}\right)
        \|\varphi\|_{H^{1/2}}
    \end{eqnarray*}
    (where we used the embedding
    $H^{1/2}\subset L^{2n/(n-1),2})$.
    From this point on, the proof proceeds exactly as in Proposition
    \ref{prop.selfadjW}
\end{proof}

\section{The limiting absorption principle}\label{sec.lap}

The essential tool in our proof will be the spectral theorem in the
following version: given a selfadjoint (unbounded) operator $A$ on
$L^{2}$ and a continuous bounded function $f(\lambda)$ on
$\mathbb{R}$, the operator $f(A)$ can be defined as
\begin{equation}\label{eq.spectral}
  f(A)\phi=-\frac1\pi\cdot L^{2}-\lim_{\epsilon\downarrow 0}
  \int f(\lambda)\Im R(\lambda+i\epsilon)\phi d\lambda
\end{equation}
for any $\phi\in L^{2}$. Here $R(z)=(A-z)^{-1}$ denotes the resolvent
operator of $A$ (see e.g. \cite{taylor}).  Under suitable assumptions on
$H$, the limit operators $R(\lambda\pm i0)=\lim_{\epsilon\downarrow 0}
R(\lambda\pm i\epsilon)$ are well defined as bounded operators in
weighted $L^{2}$ spaces; this is usually called the \emph{limiting
  absorption principle} (see below for details). Thus we have also the
simpler representation
\begin{equation}\label{eq.spectral2}
  f(A)\phi=-\frac1\pi\cdot
  \int f(\lambda)\Im R(\lambda+i0)\phi d\lambda.
\end{equation}

Recalling the definition \eqref{eq.coeff}, consider now the
operators
$$H=-\Delta+W(x,D)\equiv
-\Delta+\sum_{j=1}^{3}a_{j}(x)\partial_{j}+b(x)$$
and
$$\D_{V}=\D+V(x).$$
In Section \ref{sec.self} we proved that, under
assumptions \eqref{eq.selfadj} on $a_{j},b$ and $V(x)$, both $H$ and
$\D_{V}$ are selfadjoint operators on $L^{2}$. In particular, the
spectral formula \eqref{eq.spectral} holds for both. We shall use
the following notations: the free resolvents will be written as
$$R_{0}(z)=(-z-\Delta)^{-1},\qquad
R_{\D}(z)=(-z I_{4}+\D)^{-1}$$
while we shall use the notation $R(z)$
for both perturbed resolvents:
$$R(z)=(-z-\Delta+W)^{-1},\qquad
R(z)=(-z+\D+V)^{-1}.$$
From the context the meaning of $R(z)$ will
always be clear. Note that $R_{0}(z)$ is defined for all
$z\not\in\mathbb{R}^{+}$ while $R_{\D}(z)$ is defined for
$z\not\in\mathbb{R}$, and the same properties hold for the perturbed
resolvents.

Our first task will be to show that the stronger representation
\eqref{eq.spectral2}, i.e., the limiting absorption principle, holds
also for the perturbed operators. For $A=-\Delta$ this is a
classical result (see e.g. Agmon \cite{agmon}); here we shall use a
very precise version of the principle, due to Barcelo, Ruiz and Vega
\cite{BRV}. On the other hand, for the Dirac operator only a few
results are available, which concern the case with a nonzero mass
term (see \cite{PSU}, \cite{yamada}).

The classical results on $R_{0}$ (see \cite{agmon}) state that the
limits
\begin{equation}\label{eq.agm}
  \lim_{\epsilon\downarrow 0}R_{0}(\lambda\pm i\epsilon)=
  R_{0}(\lambda\pm i0)
\end{equation}
exist in the norm of bounded operators from $L^{2}(\xx^{s}dx)$ to
$H^{2}(\xx^{-s}dx)$ for any $s>1$; the convergence is uniform for
$\lambda$ belonging to any compact subset of $]0,+\infty[$, and the
following estimate holds
\begin{equation}\label{eq.agmon}
  \|\xx^{-s} R_{0}(\lambda\pm i0)\xx^{-s} f\|_{L^{2}}\leq
  \frac {C(s)}{\sqrt\lambda}\|f\|_{L^{2}}\qquad
  \forall\lambda>0,\ s>\frac12.
\end{equation}
In $n=3$ dimensions, the operators $R_{0}(\lambda\pm i0)$ have the
explicit representation
\begin{equation}\label{eq.expl}
  R_{0}(\lambda\pm i0)g(x)=
  \frac{1}{4\pi} \int
  \frac{e^{\pm i\sqrt\lambda |x-y|}}{|x-y|} g(y) dy,\quad\lambda\geq0.
\end{equation}
Recall also that for $\lambda<0$ we have the similar formula
\begin{equation}\label{eq.expl2a}
  R_{0}(\lambda)g(x)=
  \frac{1}{4\pi} \int
  \frac{e^{-\sqrt{|\lambda|}\; |x-y|}}{|x-y|} g(y) dy,\quad\lambda\leq0.
\end{equation}

These results were extended in \cite{BRV} to more general weights.
Introduce the norm
\begin{equation}\label{eq.brv}
  \||a(x)\||=\sup_{\mu>0}\int_{\mu}^{+\infty}
  \frac{h(r)r}{(r^{2}-\mu^{2})^{1/2}}dr\text{\ \ \ where\ \ \ }
  h(r)\equiv\sup_{|x|=r}|a(x)|.
\end{equation}
For any measurable function on $\mathbb{R}^{n}$ such that $\spt
f\subseteq\spt a$, we can consider the (semi-)norm
$$\|f\|_{L^{2}(a(x)dx)}\equiv
\|a(x)^{1/2}f\|_{L^{2}}<\infty$$
and we can define a Hilbert space
$L^{2}(a(x)dx)$ as the closure in this norm of the subspace of
$C^{\infty}_{0}$ functions with support contained in $\spt a$.  Then
we can summarize Theorems 1 and 2 in \cite{BRV} as follows:

\begin{theorem}[\cite{BRV}]\label{th.BRV}
  Let $a(x)$ be a nonnegative function on $\mathbb{R}^{n}$ with
  $\||a\||<\infty$, and denote by $R_{0}(\lambda\pm i 0)$ the limit
  operators \eqref{eq.agm}. Then the operators $R_{0}(z)$ for
  $z\not\in\mathbb{R}^{+}$ and $R_{0}(\lambda\pm i 0)$ can be extended
  to bounded operators from $L^{2}(a(x)^{-1}dx)$ to $L^{2}(a(x)dx)$,
  and the following estimates hold:
  \begin{equation}\label{eq.est1}
    \|R_{0}(\lambda\pm i 0)f\|_{L^{2}(a(x)dx)}\leq
    \frac{C}{\sqrt{|\lambda|}}\||a\||
    \cdot\|f\|_{L^{2}(a(x)^{-1}dx)},
    \qquad\lambda\neq0
  \end{equation}
  (here of course $R_{0}(\lambda\pm i 0)\equiv R_{0}(\lambda)$ for
  $\lambda<0$)
  \begin{equation}\label{eq.est2}
    \|\nabla R_{0}(\lambda\pm i 0)f\|_{L^{2}(a(x)dx)}\leq
    C\||a\||\cdot\|f\|_{L^{2}(a(x)^{-1}dx)}.
  \end{equation}
  Moreover, the limiting absorption principle holds in the weak form:
  for all $f,g\in L^{2}(a(x)^{-1}dx)$
  \begin{equation}\label{eq.lapbrv}
    \lim_{\epsilon\downarrow0} (R_{0}(\lambda\pm i \epsilon)f,g)=
    (R_{0}(\lambda\pm i 0)f,g).
  \end{equation}
\end{theorem}

\begin{remark}\label{rem.complex}
  It is not difficult to extend the estimates \eqref{eq.est1}and
  \eqref{eq.est2} to the whole complex plane. Indeed, fix two
  functions $f,g\in C^{\infty}_{0}$ with support contained in $\spt a$
  and consider on the half plane
$$S=\{z\colon\Im z>0\}$$
the holomorphic function
\begin{equation}\label{eq.Fz}
  F(z)=z^{1/2}(R_{0}(z)f,g).
\end{equation}
It is clear that $F(z)$ is continuous on $\overline{S}$ up to the
boundary, moreover it satisfies the estimate
\begin{equation}\label{eq.estFz}
  |F(x)|\leq
  C\||a\||\cdot\|f\|_{L^{2}(a(x)^{-1}dx)}
  \|g\|_{L^{2}(a(x)^{-1}dx)}
\end{equation}
on the boundary $\Im z=0$, and finally it has a polynomial growth for
$|z|\to+\infty$, as it easily follows from the explicit expression of
$R_{0}(z)$ as a convolution operator (see \cite{BRV}). By the
Phragm\'en-Lindel\"of Theorem (see e.g. \cite{Stein}) on the half
plane we immediately obtain that estimate \eqref{eq.estFz} holds on
all of $\overline{S}$. A similar argument can be applied in the lower
half plane $\Im z<0$. In conclusion we obtain
\begin{equation}\label{eq.est3}
  \|R_{0}(z)f\|_{L^{2}(a(x)dx)}\leq
  \frac{C}{\sqrt{|z|}}\;\||a\||
  \cdot\|f\|_{L^{2}(a(x)^{-1}dx)}
\end{equation}
for all $f\in L^{2}(a(x)^{-1}dx)$ (see also part (ii) in Theorem 1,
\cite{BRV}). Notice that this estimate holds on the whole complex
plane, in the sense that we apply it to $R_{0}(\lambda\pm i0)$ when
$z\in\mathbb{R}^{+}$ .

If we apply the same argument to the function
$$G(z)=(\nabla R_{0}(z)f,g)$$
we obtain in an analogous way the estimate
\begin{equation}\label{eq.est4}
  \|\nabla R_{0}(z)f\|_{L^{2}(a(x)dx)}\leq
  C\;\||a\||
  \cdot\|f\|_{L^{2}(a(x)^{-1}dx)},\qquad z\in\mathbb{C}.
\end{equation}
\end{remark}

We now specialize the theorem to a particular choice of weights.
Precisely, consider the family of functions
\begin{equation}\label{eq.wei}
  w_{\beta}(x)={|x|(|\log|x||+1)^{\beta}},\qquad \beta>1.
\end{equation}
As it is proved in \cite{BRV} (see Proposition 1), the norms
$$\||w_{\beta}^{-1}\||<+\infty$$
are finite for all $\beta >1$, hence we can apply \ref{th.BRV} with
the choice
$$a(x)=(w_{\beta}(x))^{-1}=\frac 1 {|x|(|\log|x||+1)^{\beta}}.$$
In this case it is possible to improve the above result and to obtain
a stronger version of the limiting absorption principle. To this end,
we need the following Lemma, which is inspired by \cite{agmon}:

\begin{lemma}\label{lem.agmon}
  Let $H$ be a Hilbert space, $H'$ its dual, and $H_{0}$ a second
  Hilbert space compactly embedded in $H'$. Let $T_{j},T$
  ($j=1,2,\dots$) be bounded operators in $\mathcal{L}(H,H')$ such
  that
  \begin{itemize}
  \item[(i)] $T_{j},T$ are symmetric for the pairing
    $\langle\cdot,\cdot\rangle_{H'\times H}$, i.e.,
    $$\langle Tf,g\rangle_{H'\times H}=
    \langle Tg,f\rangle_{H'\times H}\qquad \forall f,g\in H;$$
  \item[(ii)] $T_{j},T\in\mathcal{L}(H,H_{0})$ and, for some constant
    $C$ independent of $j$,
    $$\|T_{j}\|_{\mathcal{L}(H,H_{0})}\leq C.$$
\end{itemize}
Assume that
\begin{equation}\label{eq.tj}
  T_{j}f\rightharpoonup  Tf\quad
  \text{weakly in $H'$ for all $f\in H$.}
\end{equation}
Then $T_{j}\to T$ in the operator norm of $\mathcal{L}(H,H')$.
\end{lemma}

\begin{proof}
  Fix an $f\in H$; the sequence $T_{j}f$ converges weakly to $Tf$ in
  $H'$, and is bounded in $H_{0}$ by (ii), hence it admits a
  subsequence which converges in the norm of $H'$, and the limit must
  be the same i.e. $Tf$. By applying the same argument to any
  subsequence of $T_{j}f$, we conclude that the entire sequence
  $T_{j}f$ converges to $Tf$ in the norm of $H$.

  Now, let $f_{j}$ be any sequence which converges to $f$ weakly in
  $H$. Then we have for all $g\in H$
$$\langle T_{j}f_{j},g\rangle=
\langle T_{j}g,f_{j}\rangle\to \langle Tf,g\rangle$$
since $T_{j}g\to
Tg$ strongly in $H'$ and $f_{j} \rightharpoonup f$ weakly in $H$. In
other words, for any $f_{j} \rightharpoonup f$ weakly in $H$ we have
that $T_{j}f_{j} \rightharpoonup Tf$ weakly in $H'$. But, as in the
first step, we can remark that the sequence $T_{j}f_{j}$ is bounded in
$H_{0}$ and by compact embedding we obtain that the convergence is
strong: $T_{j}f_{j}\to Tf$ in the norm of $H'$.

By the same argument we obtain that, for any $f_{j} \rightharpoonup f$
weakly in $H$, the sequence $Tf_{j}$ converges to $Tf$ in the norm of
$H'$.

Finally, assume by contradiction that $T_{j}$ does not converge to $T$
in the operator norm of $\mathcal{L}(H,H')$. This means that we can
find a sequence $f_{j}\in H$ with norm $\|f_{j}\|_{H}=1$ such that
$$\|T_{j}f_{j}-Tf_{j}\|_{H'}>\epsilon>0$$
for some $\epsilon$ independent of $j$. By extracting a subsequence we
can assume that $f_{j}\rightharpoonup f$ weakly in $H$, and by the
above steps we immediately obtain a contradiction.
\end{proof}

Then we can prove:

\begin{proposition}\label{prop.lap}
  Let $w_{\beta}(x)$, $x\in\mathbb{R}^{n}$ one of the radial weights
  \eqref{eq.wei} for some fixed $\beta >1$. Then, for all
  $\lambda\neq0$, the limits
  \begin{equation}\label{eq.lapw}
    \lim_{\epsilon\downarrow 0}
    R_{0}(\lambda\pm i\epsilon)=
    R_{0}(\lambda\pm i0)
  \end{equation}
  exist in the norm of bounded operators from $L^{2}(w_{\beta}(x)dx)$
  to $H^{2}(w_{\beta}(x)^{-1}dx)$ and satisfy the estimates
  \begin{equation}\label{eq.est5}
    \|R_{0}(\lambda\pm i0)f\|_{L^{2}(w_{\beta}^{-1}dx)}\leq
    \frac{C(b)}{\sqrt{|\lambda|}}
    \;\|f\|_{L^{2}(w_{\beta}dx)},\qquad
    \forall \lambda\neq0,
  \end{equation}
  \begin{equation}\label{eq.est6}
    \|\nabla R_{0}(\lambda\pm i0)f\|_{L^{2}(w_{\beta}^{-1}dx)}\leq
    C(b)
    \;\|f\|_{L^{2}(w_{\beta}dx)}.
  \end{equation}
\end{proposition}

\begin{proof}
  We apply Lemma \ref{lem.agmon} with the choices:
  $H=L^{2}(w_{\beta}(x)dx)$, and hence $H'=L^{2}(w_{\beta}(x)^{-1}dx)$
  with the standard $L^{2}$ pairing;
  $H_{0}=H^{1}(w_{\beta_{0}}(x)^{-1}dx)$ for some arbitrary
  $\beta_{0}$ with $\beta > \beta_{0}>1$; the norm of $H_{0}$ of
  course is
$$\|f\|^{2}_{H_{0}}=\|w_{\beta_{0}}^{-1/2}f\|_{L^{2}}^{2}
+\|w_{\beta_{0}}^{-1/2}\nabla f\|_{L^{2}}^{2}.$$
Finally, as operators
$T_{j}$ we shall take (any subsequence of) the resolvent operators
$R_{0}(\lambda\pm i\epsilon)$ as $\epsilon\downarrow0$, while
$T=R_{0}(\lambda\pm i0)$, for some fixed $\lambda\in\mathbb{R}$.

We now check the assumptions of the lemma. The compact embedding of
$H_{0}$ into $H'$ is clear. Also the symmetry of the operators in the
sense of (i) is evident. The uniform bounds on $T_{j},T$ as bounded
operators from $H$ to $H'$ are simply the estimates \eqref{eq.est3},
\eqref{eq.est4} applied with the choice $a(x)=w_{\beta}(x)^{-1}$.  But
it is clear that the estimate \eqref{eq.est3} implies also the
following estimate
\begin{equation}\label{eq.est7}
  \|R_{0}(z)f\|_{L^{2}(w_{\beta_{0}}^{-1}dx)}\leq
  \frac{C(\beta_{0})}{\sqrt{|z|}}
  \;\|f\|_{L^{2}(w_{\beta}dx)},\qquad
  \forall z\neq0,
\end{equation}
which is only apparently stronger, in view of the trivial embedding
$$L^{2}(w_{\beta}dx)\subseteq
L^{2}(w_{\beta_{0}}dx).$$
In a similar way we have
\begin{equation}\label{eq.est8}
  \|\nabla R_{0}(z)f\|_{L^{2}(w_{\beta_{0}}^{-1}dx)}\leq
  {C(\beta_{0})}
  \;\|f\|_{L^{2}(w_{ \beta}dx)}.
\end{equation}
These inequalities show that assumption (ii) of the Lemma is
satisfied.  Finally, assumption \eqref{eq.tj} is nothing but the weak
limiting absorption principle of Barcelo, Ruiz, Vega (see
\eqref{eq.lapbrv}).

In conclusion, Lemma \ref{lem.agmon} implies that the limit
\eqref{eq.lapw} exists in the norm of bounded operators from
$L^{2}(w_{\beta}dx)$ to $L^{2}(w_{\beta}^{-1}dx)$. Moreover, by the
identity
$$\Delta R_{0}(z)=-I-z R_{0}(z)$$
we obtain that the limit exists also in the norm of bounded operators
from $L^{2}(w_{\beta}dx)$ to $H^{2}(w_{\beta}^{-1}dx)$. The estimates
\eqref{eq.est5} and \eqref{eq.est6} follow from the corresponding
estimates for general $z$.
\end{proof}

\subsection{The limiting absorption principle for the magnetic
  laplacian}\label{ssec.magn}

In what follows, we shall focus on the case $n=3$ exclusively.  We
follow the standard approach, based on the resolvent identity
\begin{equation*}
  R(z)=(-z-\Delta+W(x,D))^{-1}=R_{0}(z)(I+WR_{0}(z))^{-1}.
\end{equation*}
Thus the main step of the proof will consist in inverting the
operator $I+WR_{0}$ in suitable weighted spaces. We shall assume
that the coefficients $a_{j}(x)$ and $b(x)$ in $W(x,D)$, defined as
in \eqref{eq.coeffi}, satisfy the assumptions
\begin{equation}\label{eq.assab}
  |a_{j}(x)|\leq\frac{C_{0}}{|x|\xx^{s}
    (|\log |x\|+1)^{\beta}},\qquad
  |b(x)|\leq\frac {C_{0}} {|x|^{2}(|\log |x\|+1)^{\beta}}
\end{equation}
for some $s\in[0,1]$, $\beta>1$ and some constant $C_{0}$ small
enough.

Our result is the following:

\begin{proposition}\label{prop.lapW}
  Assume the coefficients of $W(x,D)=\sum a_{j}(x)\partial_{j}+b(x)$
  satisfy \eqref{eq.coeffi}
   \eqref{eq.assab} for some $C_{0}$ small enough, some
  $s\in[0,1]$ and some $\beta>1$.

  Then the operator $I+WR_{0}$ is invertible on the
  weighted space
  $L^{2}(w_{\beta}(x)\xx^{2s}dx)$, and the inverse operators
  $(I+WR_{0}(z))^{-1}$ are uniformly bounded for all $z\in\mathbb{C}$.
  Moreover, the strong limiting absorption principle holds for $R(z)$,
  in the following sense:
  \begin{itemize}
  \item[(i)] the boundary values
    \begin{equation}\label{eq.agmW}
      \lim_{\epsilon\downarrow 0}
      R(\lambda\pm i\epsilon)=
      R(\lambda\pm i0)
    \end{equation}
    exist in the norm of bounded operators from
    $L^{2}(w_{\beta}(x)dx)$ to $H^{2}(w_{\beta}^{-1}(x)dx)$;
  \item[(ii)] the following estimate
    \begin{equation}\label{eq.estR}
      \|R(z)f\|_{L^{2}(w_{\beta}(x)dx)}\leq
      \frac{C(\beta)}{\sqrt{|z|}}
      \cdot\|f\|_{L^{2}(w_{\beta}(x)^{-1}dx)}
    \end{equation}
    holds for all $z\in\mathbb{C}$, $z\neq0$.
  \end{itemize}
\end{proposition}

\begin{remark}\label{rem.s}
  In the case $s=0$ we recover exactly the strong limiting absorption
  principle proved in Proposition \ref{prop.lap} above for the free
  operator $R_{0}$. The additional weight $\xx^{s}$ was considered in
  view of the estimates that will be needed in the following section.
\end{remark}

\begin{proof}
  Consider the operator
  \begin{equation*}
    W(x,D)R_{0}(z)f=
    \sum a_{j}(x)\partial_{j}R_{0}(z)f+
    b(x)R_{0}(z)f;
  \end{equation*}
  we estimate the two terms separately.

  First of all we have
  \begin{equation*}
    \|w_{\beta}^{1/2}\xx^{s}a_{j}(x)
    \partial_{j}R_{0}f\|_{L^{2}}\leq
    \|w_{\beta}\xx^{s}a_{j}\|_{L^{\infty}}
    \|w_{\beta}^{-1/2}\partial_{j}R_{0}f\|_{L^{2}}\leq
    C_{0}\|w_{\beta}^{1/2}f\|_{L^{2}}
  \end{equation*}
  by estimate \eqref{eq.est8}, and this implies trivially
  \begin{equation}\label{eq.wR0}
    \|w_{\beta}^{1/2}\xx^{s}a_{j}(x)
    \partial_{j}R_{0}f\|_{L^{2}}\leq
    C_{0}\|w_{\beta}^{1/2}\xx^{s}f\|_{L^{2}}.
  \end{equation}

  In order to estimate the electric term, we recall that, from the
  explicit expression of the free resolvent, we can write
  \begin{equation*}
    |R_{0}(z)f|\leq \frac1{4\pi} \left|\frac1{|x|}*|f|\right|.
  \end{equation*}

  Then we have
  \begin{equation}\label{eq.wbR0}
    \|w_{\beta}^{1/2}b(x)R_{0}(z)f\|_{L^{2}}\leq
    \|w_{\beta}^{1/2}b(x)\|_{L^{2}}
    \|R_{0}(z)f\|_{L^{\infty}}\leq
    \|w_{\beta}^{1/2}b(x)\|_{L^{2}}\cdot
    C\left\|\frac1{|x|}*|f|\right\|_{L^{\infty}}.
  \end{equation}
  Recalling Young and H\"older inequalities in Lorentz spaces
  (see Theorems \ref{th.holderlorentz}, \ref{th.younglorentz}),
  we have
  \begin{equation*}
    \left\|\frac1{|x|}*|f|\right\|_{L^{\infty}}\leq
    C\|f\|_{L^{3/2,1}}=
    C\|w_{\beta}^{-1/2}w_{\beta}^{1/2}f\|_{L^{3/2,1}}\leq
    C\|w_{\beta}^{-1/2}\|_{L^{6,2}}
    \| w_{\beta}^{1/2}f\|_{L^{2}}.
  \end{equation*}
  Since $w_{\beta}^{-1/2}\in L^{6,2}$ for any $\beta>1$
  (Proposition \ref{prop.wbeta}), \eqref{eq.wbR0} gives
  \begin{equation*}
    \|w_{\beta}^{1/2}b(x)R_{0}(z)f\|_{L^{2}}\leq
    C\|w_{\beta}^{1/2}b(x)\|_{L^{2}}\cdot
    \| w_{\beta}^{1/2}f\|_{L^{2}}.
  \end{equation*}
  Now, by assumption \eqref{eq.assab} on $b(x)$ we have easily
  \begin{equation*}
    \|w_{\beta}^{1/2}b(x)\|_{L^{2}}\leq C C_{0}
  \end{equation*}
  and we conclude that
  \begin{equation}\label{eq.wbR0fin}
    \|w_{\beta}^{1/2}b(x)R_{0}(z)f\|_{L^{2}}\leq
    CC_{0}\cdot
    \| w_{\beta}^{1/2}f\|_{L^{2}}.
  \end{equation}

  In a similar way we have
  \begin{equation}\label{eq.wbR1}
    \|w_{\beta}^{1/2}\xx b R_{0}(z)f\|_{L^{2}}\leq
    \|w_{\beta}^{1/2}\xx b\|_{L^{6}}
    \|R_{0}(z)f\|_{L^{3}}\leq
    \|w_{\beta}^{1/2}\xx b\|_{L^{6}}\cdot
    C\left\|\frac1{|x|}*|f|\right\|_{L^{3}}
  \end{equation}
  and
  \begin{align*}
    \left\|\frac1{|x|}*|f|\right\|_{L^{3}}\leq
    C\|f\|_{L^{1}}=&
    C\|w_{\beta}^{-1/2}\xx^{-1}
    w_{\beta}^{1/2}\xx f\|_{L^{1}}\\
    \leq&
    C\|w_{\beta}^{-1/2}\xx^{-1}\|_{L^{2}}
    \| w_{\beta}^{1/2}\xx f\|_{L^{2}}.
  \end{align*}
  As above, we notice that $w_{\beta}^{-1/2}\xx^{-1}\in L^{2}$ for any
  $\beta>1$, hence we have from \eqref{eq.wbR1}
  \begin{equation*}
    \|w_{\beta}^{1/2}\xx b R_{0}(z)f\|_{L^{2}}\leq
    C\|w_{\beta}^{1/2}\xx b\|_{L^{6}}\cdot
    \|w_{\beta}^{1/2}\xx f\|_{L^{2}}.
  \end{equation*}
  Assumption \eqref{eq.assab} guarantees that
  \begin{equation*}
    \|w_{\beta}^{1/2}\xx b(x)\|_{L^{6}}\leq C C_{0}
  \end{equation*}
  and, in conclusion,
  \begin{equation}\label{eq.wbR2}
    \|\xx w_{\beta}^{1/2} b(x) R_{0}(z)f\|_{L^{2}}\leq
    CC_{0}\cdot\|\xx w_{\beta}^{1/2} f\|_{L^{2}}
  \end{equation}
  If we interpolate between \eqref{eq.wbR0fin} and \eqref{eq.wbR2}, we
  obtain the estimate
  \begin{equation}\label{eq.wbRs}
    \|\xx^{s} w_{\beta}^{1/2} b(x) R_{0}(z)f\|_{L^{2}}\leq
    CC_{0}\cdot\|\xx^{s} w_{\beta}^{1/2} f\|_{L^{2}}
  \end{equation}

  Summing up, from estimates \eqref{eq.wR0} and \eqref{eq.wbRs} we get
  for all $z\in\mathbb{C}$
  \begin{equation}\label{eq.estfin}
    \|\xx^{s} w_{\beta}^{1/2} WR_{0}(z)f\|_{L^{2}}\leq
    CC_{0}\cdot\|\xx^{s} w_{\beta}^{1/2} f\|_{L^{2}}.
  \end{equation}
  Then it is clear that we can invert the operator $I+WR_{0}$ by a
  Neumann series on the space $L^{2}(\xx^{2s} w_{\beta}dx)$.  Hence,
  the standard representation
  \begin{equation}\label{eq.LS}
    R(z)=R_{0}(z)(I+WR_{0}(z))^{-1}
  \end{equation}
  is valid. To conclude the proof of the Proposition, it is now
  sufficient to remark that, from property \eqref{eq.lapw} of
  Proposition \ref{prop.lap} and the uniform bounds on the norm of
  $(I+WR_{0}(z))^{-1}$ we have just obtained (for $s=0$), the limits
  in \ref{eq.agmW} exist in a weak sense. Proceeding as in the proof
  of Proposition \ref{prop.lap}, using Lemma \ref{lem.agmon}, we
  deduce (i). Finally, (ii) is a consequence of \eqref{eq.LS} and the
  corresponding estimate \eqref{eq.est7} for $R_{0}$.
\end{proof}

\begin{remark}\label{rem.AB}
Note that the assumptions of the preceding proposition can be
expressed in terms of the original coefficients $A,B$ as
follows:
\begin{equation}\label{eq.assabb}
  |A(x)|\leq\frac{C_{0}}{|x|\xx^{s}
    (|\log |x\|+1)^{\beta}},\qquad
  |\nabla A(x)|+
  |B(x)|\leq\frac {C_{0}} {|x|^{2}(|\log |x\|+1)^{\beta}}
\end{equation}
for some $\beta>1$ and a constant $C_{0}>0$ small enough.
\end{remark}

\subsection{The limiting absorption principle for the Dirac operator
  and its perturbation}
\label{ssec.dirac}

In this section we will study the limiting absorption principle for
the massless Dirac operator $\D$; this property was studied by Yamada
in \cite{yamada} for the operator with mass. Moreover, as in the case
of the magnetic Laplacian, we will extend this result to the perturbed
operator $\D_V=\D+V(x)$, under a suitable assumption on the potential
$V$.

It is well known that the spectrum of the free operator $\D$ is the
whole real line. Due to the relation $\D^2=-\Delta I_4$, we
immediately obtain the representation
\begin{equation}\label{eq.diracres}
  R_{\D}(z)=R_0(z^2)(\D+zI_4),
\end{equation}
for all $z\in\mathcal{C}$ with $\Re z=0$. Using this formula and the
Proposition \ref{prop.lap}, we easily prove the following:

\begin{proposition}\label{prop.RD}
  Let $w_{\beta}(x)$, $x\in\mathbb{R}^{3}$ be defined as in
  \eqref{eq.wei}, for some fixed $\beta >1$.  Then, for all
  $\lambda\in\mathbb{R}$, the limits
  \begin{equation}\label{eq.lapdirac}
    \lim_{\epsilon\downarrow 0}
    R_{\D}(\lambda\pm i\epsilon)=
    R_{\D}(\lambda\pm i0):=R_0(\lambda^2\pm i0)(\D+\lambda I_4)
  \end{equation}
  exist in the norm of bounded operators from $L^{2}(w_{\beta}(x)dx)$
  to $H^1(w_{\beta}(x)^{-1}dx)$ and satisfy the estimate
  \begin{equation}\label{eq.estD}
    \|R_{\D}(z)f\|_{L^2(w_{\beta}(x)^{-1}dx)}\leq\|f\|_{L^2(w_{\beta}(x)dx)},
  \end{equation}
  for all $z\in\mathbb{C}$.  Moreover, we have the explicit
  representation
  \begin{eqnarray}
    R_{\D}(\lambda\pm i0)f & = & \frac{i|\lambda|}{4\pi}
    \int_{\R^3}\frac{e^{i|\lambda|\cdot|x-y|}}{|x-y|}\left(
      I_4-\sum_{j=1}^{3}\alpha_j\frac{x_j-y_j}{|x-y|}\right)f(y)\diff
    y\nonumber
    \\ \  & \  &
    +\frac{1}{4\pi}\int_{\R^3}
     \frac{e^{i|\lambda|\cdot|x-y|}}{|x-y|^2}\sum_{j=1}^{3}
    \alpha_j\frac{x_j-y_j}{|x-y|}f(y)\diff y.\label{eq.diracres2}.
  \end{eqnarray}
\end{proposition}

\begin{proof}
  The strong convergence of $R_{\D}(\lambda\pm i\epsilon)$ to
  $R_{\D}(\lambda\pm i0)$ in the space of bounded operators from
  $L^{2}(w_{\beta}(x)dx)$ to $H^1(w_{\beta}(x)^{-1}dx)$ is obtained by
  interpolation using the property \eqref{eq.lapw} and the
  representation \eqref{eq.diracres}; estimate \eqref{eq.estD}
  immediately follows from \eqref{eq.diracres} and the estimates
  \eqref{eq.est5}, \eqref{eq.est6}, \eqref{eq.est7}, \eqref{eq.est8}.
  In conclusion, recalling the explicit representation \eqref{eq.expl}
  for $R_0(\lambda\pm i0)$, after an integration by parts we get the
  formula \eqref{eq.diracres2} and this concludes the proof.
\end{proof}

At this point, we will proceed in a similar way to the case of the
perturbed Laplacian and we will prove that it is possible to extend
the above result to small electric perturbations of the free Dirac
operator. As for the magnetic coefficients of $W(x,D)$, we need to
assume that the potential $V$ satisfies
\begin{equation}\label{eq.assV}
  |V(x)|\leq\frac{C_{0}}{|x|\xx^{s}
    (|\log |x\|+1)^{\beta}},
\end{equation}
for some $s\in[0,1]$, $\beta>1$ and some constant $C_{0}$ small
enough. We prove the following result:

\begin{proposition}\label{prop.lapRV}
  Assume the potential $V$ satisfies \eqref{eq.assV} for some $C_{0}$
  sufficiently small, some $s\in[0,1]$ and some $\beta>1$.

  Then the operator $I+VR_{\D}$ is invertible on the weighted space
  $L^{2}(w_{\beta}(x)\xx^{2s}dx)$, and the inverse operators
  $(I+VR_{\D}(z))^{-1}$ are uniformly bounded for all
  $z\in\mathbb{C}$. Moreover, the strong limiting absorption principle
  holds for $R(z)$, in the following sense:
  \begin{itemize}
  \item[(i)] the limits
    \begin{equation}\label{eq.agmV}
      \lim_{\epsilon\downarrow 0}
      R(\lambda\pm i\epsilon)=
      R(\lambda\pm i0)
    \end{equation}
    exist in the norm of bounded operators from
    $L^{2}(w_{\beta}(x)dx)$ to $H^{1}(w_{\beta}^{-1}(x)dx)$;
  \item[(ii)] the following estimate
    \begin{equation}\label{eq.estRV}
      \|R(z)f\|_{L^{2}(w_{\beta}(x)^{-1}dx)}\leq
      C(\beta)
      \cdot\|f\|_{L^{2}(w_{\beta}(x)dx)}
    \end{equation}
    holds for all $z\in\mathbb{C}$, $z\neq0$.
  \end{itemize}
\end{proposition}

\begin{proof}
  The argument is the same of the proof of Proposition \ref{prop.lapW}
  for the magnetic part of $W$. First we observe that, by hypothesis
  \eqref{eq.assV}, we have
  \begin{equation*}
    \|w_\beta^{1/2}\xx^s V(x)R_{\D}f\|_{L^2}\leq\|w_\beta\xx^s
    V(x)\|_{L^\infty}\|w_\beta^{-1/2}R_{\D}f\|_{L^2}\leq C_0\cdot
    \|w_\beta^{-1/2}f\|_{L^2}.
  \end{equation*}
  Hence we obtain the estimate
  \begin{equation*}
    \|w_\beta^{1/2}\xx^s V(x)R_{\D}(z)f\|_{L^2}\leq
    \|w_\beta^{1/2}\xx^s f\|_{L^2},
  \end{equation*}
  uniformly in $z\in\mathbb{C}$; thus we can invert the operator
  $I+VR_{\D}$ by a Neumann series on the space $L^2(w_\beta dx)$.
  Again, we can exploit the representation
  \begin{equation}\label{eq.LSdirac}
    R(z)=R_{\D}(z)(I+VR_{\D}(z))^{-1}.
  \end{equation}
  By property \eqref{eq.lapdirac} of Proposition \ref{prop.RD} and the
  uniform bounds of $(I+VR_\D)^{-1}$, it follows that the limits in
  \eqref{eq.agmV} exist in a weak sense. Then we can procede as in the
  previous cases, using Lemma \ref{lem.agmon} and obtain (i). In
  conclusion, the estimate (ii) is an immediate consequence of
  \eqref{eq.LSdirac} and the inequality \eqref{eq.estD}. This
  concludes the proof.
\end{proof}

In the following we shall also need a weaker version of the last
result: we shall require that $V$ satisfies
\begin{equation}\label{eq.assV2}
  |V(x)|\leq\frac{C_{0}}{|x|^{1/2}\xx^{s}
    (|\log|x\|+1)^{\beta/2}},
\end{equation}
for some $s>\frac12$, $\beta>1$ and some constant $C_{0}$ small
enough. Then we have

\begin{corollary}\label{cor.lapRV2}
  Assume the potential $V$ satisfies \eqref{eq.assab} for some $C_{0}$
  sufficiently small, $s>\frac12$ and $\beta>1$.

  Then the operators $I+VR_{\D}$ are invertible on the space
  $L^{2}(\xx^{2s}dx),$ and the inverse operators $(I+VR_{\D}(z))^{-1}$
  are uniformly bounded for all $z\in\mathbb{C}$. Moreover, the strong
  limiting absorption principle holds for $R(z)$, in the following
  sense:
  \begin{itemize}
  \item[(i)] the limits
    \begin{equation}\label{eq.agmV2}
      \lim_{\epsilon\downarrow 0}
      R(\lambda\pm i\epsilon)=
      R(\lambda\pm i0)
    \end{equation}
    exist in the norm of bounded operators from $L^{2}(\xx^{2s}dx)$ to
    $ H^{1}(\xx^{-2s}dx);$
  \item[(ii)] the following estimate
    \begin{equation}\label{eq.estRV2}
      \|R(z)f\|_{L^{2}(\xx^{-2s}dx)}\leq
      C\cdot\|f\|_{L^{2}(\xx^{2s}dx)}
    \end{equation}
    holds for all $z\in\mathbb{C}$, $z\neq0$.
  \end{itemize}
\end{corollary}

\begin{proof}
  The proof is analogous to the proof of Proposition
  \ref{prop.lapRV}. Indeed, from estimate \eqref{eq.estD} and
  assumption \eqref{eq.assV2} we have immediately
  \begin{equation*}
    \|\xx^{s}VR_{\D}\|_{L^{2}}\leq
    \|\xx^{s}w_{\beta}^{1/2}V\|_{L^{\infty}}
    \|w_{\beta}^{-1/2}R_{\D}f\|_{L^{2}}\leq
    C_{0}\|w_{\beta}^{1/2}f\|_{L^{2}}
  \end{equation*}
  and by the trivial inequality
  \begin{equation*}
    w_{\beta}^{1/2}\leq  C_{s}\xx^{s},
  \end{equation*}
  valid for all $s>1/2$, we conclude that
  \begin{equation*}
    \|\xx^{s}VR_{\D}\|_{L^{2}}\leq
    C_{0} \|\xx^{s}f\|_{L^{2}}.
  \end{equation*}
  Thus we can again invert $(I+VR_{\D})$ with a Neumann series, and
  proceeding exactly as before we obtain the proof of the Corollary.
\end{proof}

\section{Resolvent Estimates}\label{sec.res.est}

In this section we prepare the crucial resolvent estimates that will
be used in the proof of the main results. In order to use the
spectral formula, we need estimates on the perturbed resolvent
operators and their derivatives with respect to $\lambda$ as bounded
operators from suitable weighted $L^p$ spaces to $L^\infty$. We
shall use the H\"older and Young inequalities in Lorentz spaces
extensively; for the convenience of the reader, we give a sketch of
the main usefule results in the Appendix \ref{appendix}.

We consider first the resolvent of the magnetic laplacian. We recall
that, by Proposition \ref{prop.lapW}, the operators $R(\lambda\pm
i0)=R_0(\lambda\pm i0)(I+W(x,D)R_0(\lambda\pm i0))^{-1}$ are well
defined as bounded operators from $L^2(w_\beta(x)dx)$ to
$H^2(w_\beta(x)^{-1}dx)$; moreover, we have the explicit
representation \eqref{eq.expl}. Our first result is the following:

\begin{lemma}\label{lem.stimeW}
  Let $R(\lambda\pm i0)=R_0(\lambda\pm i0)(I+W(x,D)R_0(\lambda\pm
  i0))^{-1}$ be the resolvent of $-\Delta+W$ and assume the
  coefficients of $W(x,D)=\sum a_j(x)\partial_j+b(x)$ satisfy
  \eqref{eq.assab}. Then, for all $\lambda\geq0$, the following
  estimates hold:
  \begin{equation}\label{eq.stimaRW}
    \|R(\lambda\pm i0)f\|_{L^\infty}\leq
    C\|w_\beta^{1/2}f\|_{L^2},
  \end{equation}
  \begin{equation}\label{eq.stimaRW2}
    \|\partial_\lambda R(\lambda\pm i0)f\|_{L^\infty}
    \leq C\left(1+\frac1{\sqrt\lambda}\right)
    \|\xx w_\beta^{1/2}f\|_{L^2}.
  \end{equation}
\end{lemma}
\begin{proof}
  The estimate \eqref{eq.stimaRW} is the easiest one. In fact, by
  formula \eqref{eq.LS} and the explicit representation
  \eqref{eq.expl} for $R_0$, we obtain
  \begin{equation*}
    \|R(\lambda\pm i0)f\|_{L^\infty}\leq
    C\cdot\|\frac{1}{|x|}*|(I+WR_0)^{-1}f|\|_{L^\infty};
  \end{equation*}
  using Young inequality in Lorentz spaces, we get
  \begin{eqnarray*}
    \|R(\lambda\pm i0)f\|_{L^\infty} & \leq &
    \|(I+WR_0)^{-1}f\|_{L^{3/2,1}} \\ \  & \leq &
    \|w_\beta(x)^{-1/2}w_\beta(x)^{1/2}(I+WR_0)^{-1}f\|_{L^{3/2,1}}
    \\ \  & \leq &
    \|w_\beta(x)^{-1/2}\|_{L^{6,2}}\|w_\beta(x)^{1/2}(I+WR_0)^{-1}f\|_{L^2}.
  \end{eqnarray*}
  The uniform bound for the operators $(I+WR_0)^{-1}$ proved in
  Proposition \ref{prop.lapW} and the observation that
  $w_\beta^{-1/2}\in L^{6,2}$, for all $\beta>1$ (see Proposition \ref{prop.wbeta})
  are sufficient now to conclude the proof of estimate
  \eqref{eq.stimaRW}.

  In order to proceed with the proof of \eqref{eq.stimaRW2} we observe
  that from \eqref{eq.expl} we immediately obtain the following
  explicit representations, for all $\lambda>0$:
  \begin{equation}\label{eq.expl2}
    \partial_\lambda R_0(\lambda\pm i0)f=R_0^2(\lambda\pm i0)f
    =\pm\frac i{8\pi\sqrt\lambda}\int_0^\infty
    e^{\pm i\sqrt\lambda|x-y|}
    f(y) dy,
  \end{equation}
  \begin{equation}\label{eq.expl3}
    \partial_j R_0^2(\lambda\pm
    i0)f=\pm\frac1{8\pi}\int_0^\infty e^{\pm i\sqrt\lambda|x-y|}
    \sum\frac{x_j-y_j}{|x-y|}f(y) dy.
  \end{equation}
  At this point, differentiating in \eqref{eq.LS} we get
  \begin{equation}\label{eq.explV}
    \partial_\lambda R(\lambda\pm i0) =A+B
  \end{equation}
  where
  \begin{equation*}
    A=R_0^2(\lambda\pm i0)
           (I+WR_0(\lambda\pm i0))^{-1}
\end{equation*}
and
\begin{equation*}
    B= R_0(\lambda\pm i0)
           (I+WR_0(\lambda\pm i0))^{-1}
           WR_0^2(\lambda\pm i0)
          (I+WR_0(\lambda\pm i0))^{-1}.
\end{equation*}
  We treat separately the two terms. By \eqref{eq.expl2}, we estimate
  \begin{eqnarray*}
    \|Af\|_{L^\infty} & \leq & \frac
    C{\sqrt\lambda}\|(I+WR_0)^{-1}f\|_{L^1} \\ \  & \leq &
    \frac C{\sqrt\lambda}\|\xx^{-1}w_\beta(x)^{-1/2}\|_{L^2}
    \|\xx w_\beta(x)^{1/2}(I+WR_0)^{-1}f\|_{L^2}.
  \end{eqnarray*}
  We observe (Proposition \ref{prop.wbeta}) that $\xx^{-1}w_\beta(x)^{-1/2}\in L^2$
  for all $\beta>1$ and, by the uniform bound for the norms of
  $(I+WR_0)^{-1}$ in the space of bounded operators onto $L^2(\xx
  w_{\beta}(x)dx)$ for (see Proposition \ref{prop.lapW}), we conclude
  that, for some $C>0$
  \begin{equation}\label{eq.stimaA}
    \|Af\|_{L^\infty}\leq\frac C{\sqrt\lambda}\|\xx
    w_\beta(x)^{1/2}f\|_{L^2}.
  \end{equation}
  For the estimate of the term $B$, we start with some computation on
  the operator $WR_0^2$. Using the representation \eqref{eq.expl3}, we
  obtain
  \begin{equation*}
    \|w_\beta^{1/2}a_j\partial_jR_0^2f\|_{L^2}\leq
    \|w_\beta^{1/2}a_j\|_{L^2}\|\partial_jR_0^2f\|_{L^\infty}\leq
    C\cdot\|w_\beta^{1/2}a_j\|_{L^2}\|f\|_{L^1}.
  \end{equation*}
  By the above observation that
    $$
    \|f\|_{L^1}\leq\|\xx w_\beta^{1/2}(x)f\|_{L^2},
    $$
    it turns out that, if $w_\beta^{1/2}a_j\in L^2$, then
    \begin{equation}\label{eq.WR1}
      \|w_\beta(x)^{1/2}a_j(x)\partial_jR_0^2f\|_{L^2}\leq
      C\cdot\|\xx w_\beta^{1/2}(x)f\|_{L^2}.
    \end{equation}
    In a similar way, using \eqref{eq.expl2}, we have
    \begin{equation*}
      \|w_\beta^{1/2}bR_0^2f\|_{L^2}\leq
      \|w_\beta^{1/2}b(x)\|_{L^2}\|R_0^2f\|_{L^\infty}\leq
      \frac C{\sqrt\lambda}\|w_\beta^{1/2}b\|_{L^2}\|f\|_{L^1}.
    \end{equation*}
    If we assume that $w_\beta^{1/2}b\in L^2$, we conclude that
    \begin{equation}\label{eq.WR2}
      \|w_\beta(x)^{1/2}b(x)R_0^2f\|_{L^2}\leq
      \frac C{\sqrt\lambda}\cdot\|\xx w_\beta^{1/2}(x)f\|_{L^2}.
    \end{equation}
    Inequalities \eqref{eq.WR1} and \eqref{eq.WR2} can be unified now,
    to show that, under the assumptions
    \begin{equation}\label{eq.assabdeb}
      w_\beta^{1/2}a_j\in L^2, \qquad w_\beta^{1/2}b\in L^2,
    \end{equation}
    the estimate
    \begin{equation}\label{eq.stimaWR2}
      \|w_\beta(x)^{1/2}W(x,D)R_0^2(\lambda\pm i0)f\|_{L^2}
      \leq C\left(1+\frac1{\sqrt\lambda}\right)
      \|\xx w_\beta(x)^{1/2}f\|_{L^2}
    \end{equation}
    holds, for some $C>0$. Observe that assumptions
    \eqref{eq.assabdeb} are weaker than \eqref{eq.assab}, so that they
    are obviously satisfied by the hypothesis of the Lemma.

    Now we are ready for the estimate of the term $B$. First, we use
    the representation \eqref{eq.expl} for $R_0$ to obtain
    \begin{eqnarray*}
      \|Bf\|_{L^\infty} & \leq &
      \left\|\frac{1}{|x|}*|(I+WR_0)^{-1}WR_0^2
        (I+WR_0)^{-1}f|\right\|_{L^\infty} \\ \  & \leq &
      \|(I+WR_0)^{-1}WR_0^2(I+WR_0)^{-1}f\|_{L^{3/2,1}}=:\|Tf\|_{L^{3/2,1}}.
    \end{eqnarray*}
    As before, we use the properties of the weights $w_\beta(x)$ to
    observe that
    $$
    \|g\|_{L^{3/2,1}}\leq\|w_\beta(x)^{1/2}g\|_{L^2}.
    $$
    Then, the last series of inequalities gives
    \begin{equation*}
      \|Bf\|_{L^\infty}\leq\|w_\beta(x)^{1/2}Tf\|_{L^2}.
    \end{equation*}
    Now we use the uniform bounds for the inverse operators
    $(I+WR_0)^{-1}$ (see Proposition \ref{prop.lapW}) to proceed with
    \begin{equation*}
      \|Bf\|_{L^\infty}\leq\|w_\beta(x)^{1/2}WR_0^2(I+WR_0)^{-1}f\|_{L^2};
    \end{equation*}
    finally, by inequality \eqref{eq.stimaWR2} and the above mentioned
    estimates on the norms of $(I+WR_0)^{-1}$ in the space of bounded
    operators onto $L^2(\xx w_\beta(x)^{1/2}dx)$, we obtain the
    estimate
    \begin{equation}\label{eq.stimaB}
      \|Bf\|_{L^\infty}\leq C\left(1+\frac1{\sqrt\lambda}\right)
      \|\xx w_\beta(x)^{1/2}f\|_{L^2}.
    \end{equation}
    In conclusion, estimates \eqref{eq.stimaA}, \eqref{eq.stimaB} and
    the representation \eqref{eq.explV} conclude the proof of
    \eqref{eq.stimaRW2} and the Lemma.
  \end{proof}

\begin{remark}\label{rem.spec1}
  The limiting absorption principle allows us to rewrite the spectral
  formula in the following way: for any (smooth, compactly supported)
  function $\phi(\lambda)$ on $\mathbb{R}$, and any test function $f$,
  \begin{equation}\label{eq.spec1}
    \phi(-\Delta+W)f=
    \int_{0}^{+\infty}\phi(\lambda)
    \Im R(\lambda+i0)fd\lambda.
  \end{equation}
  where the integral is restricted to the positive real axis since of
  course $\Im R(\lambda)=0$ for negative $\lambda$.

  The resolvent estimates just proved imply that we can integrate by
  parts in the above formula, i.e., if
  \begin{equation*}
    \phi({\lambda})=\psi'(\lambda)
  \end{equation*}
  then
  \begin{align}\label{eq.intp}
    \phi(-\Delta+W)f=& \int_{0}^{+\infty}\psi'(\lambda)
    \Im R(\lambda+i0)fd\lambda\\
    =&- \int_{0}^{+\infty}\psi(\lambda) \partial_{\lambda}\Im
    R(\lambda+i0)fd\lambda \nonumber
  \end{align}
  The problems arising from the singularity at $\lambda=0$ are easily
  overcome.  To prove this, consider a cutoff function $\chi(\lambda)$
  supported in $[-L,L]$, and write
  \begin{equation*}
    \phi(-\Delta+W)f=\lim_{L\to+\infty}
    \int_{0}^{+\infty}\phi(\lambda)(1-\chi(\lambda L))
    \Im R(\lambda+i0)fd\lambda
  \end{equation*}
  whence
  \begin{align*}
    \phi(-\Delta+W)f=&-\lim_{L\to{+\infty}} L
    \int_{0}^{1/L}\psi(\lambda)\chi'(\lambda L)
    \Im R(\lambda+i0)fd\lambda\\
    &- \lim_{L\to{+\infty}} \int_{0}^{+\infty}(1-\chi(\lambda L))
    \psi(\lambda) \partial_{\lambda}
    \Im R(\lambda+i0)fd\lambda.\\
    =&u_{L}+v_{L}.
  \end{align*}
  The last term $v_{L}$ converges to \eqref{eq.intp} uniformly, thanks
  to estimate \eqref{eq.stimaRW2} (and Lebesgue's dominated
  convergence theorem), hence it is clear that
  $u_{L}=\phi(-\Delta+W)f-v_{L}$ also converges uniformly, and it will
  be sufficient to show that its limit is 0, e.g., in distribution
  sense. To estimate the integral
  \begin{equation*}
    u_{L}=-L
    \int_{0}^{1/L}\psi(\lambda)\chi'(\lambda L)
    \Im R(\lambda+i0)fd\lambda
  \end{equation*}
  we can use the identity
  \begin{equation}\label{eq.imR}
    \Im R(\lambda+i0)=(I+R_{0}(\lambda-i0)W)^{-1}
    \Im R_{0}(\lambda+i0)
    (I+WR_{0}(\lambda+i0))^{-1}.
  \end{equation}
  Consider then the $L^{2}$ product
  \begin{equation*}
    (\Im R(\lambda+i0)f,g)=
    (\Im R_{0}(\lambda+i0)
    (I+WR_{0}(\lambda+i0))^{-1}f,(I+WR_{0}(\lambda+i0))^{-1}g).
  \end{equation*}
  From the explicit formula
  \begin{equation*}
    \Im R_{0}(\lambda+i0)h=C\int
    \frac{\sin(\sqrt\lambda|x-y|)}{|x-y|}h(y)dy
  \end{equation*}
  we have
  \begin{equation*}
    |\Im R_{0}(\lambda+i0)h|\leq C\sqrt\lambda\int
    |h(y)|dy\
  \end{equation*}
  which implies
  \begin{equation*}
    \|\Im R_{0}(\lambda+i0)h\|_{L^{\infty}}
    \leq C\sqrt\lambda\|h\|_{L^{1}}\leq
    C\sqrt\lambda\|\xx w_{\beta}^{1/2}h\|_{L^{2}}
  \end{equation*}
  for any $\beta>1$.  Recalling now the uniform bound for
  $(I+WR_{0}(\lambda+i0))^{-1}$ in Proposition \ref{prop.lapW} in the
  weighted $L^{2}$ norms with weight $\xx w_{\beta}^{1/2}$, we obtain
  easily
  \begin{equation*}
    |(\Im R(\lambda+i0)f,g)|\leq C
    \sqrt\lambda
    \|\xx w_{\beta}^{1/2}f\|_{L^{2}}
    \|\xx w_{\beta}^{1/2}g\|_{L^{2}}.
  \end{equation*}
  From this estimate it is easy to prove that
  \begin{equation*}
    (u_{L},g)=-L\int_{0}^{1/L}\psi(\lambda)\chi'(\lambda L)
    (\Im R(\lambda+i0)f,g)d\lambda\to0
  \end{equation*}
  as $L\to+\infty$, which concludes the argument.
\end{remark}

We will prove now an analogue of Lemma \ref{lem.stimeW} for the Dirac
operator. In what follows, $R(z)=(-zI_4+\D+V)^{-1}$ denotes the
resolvent of the perturbed Dirac operator. Our approach here will be
slightly different: we shall use the formula
\begin{equation}\label{eq.LS3}
  R(z)=R_\D(z)+R_\D(z)V(x)R_\D(z)(I+V(x)R_\D(z))^{-1},
\end{equation}
valid for all $z\in\mathbb{C}$ (to be interpreted of course, for
$z=\lambda\in\mathbb{R}$, as the extended resolvents
$R(\lambda):=R(\lambda\pm i0)$ on the weighted $L^2$ spaces, as given
by Proposition \ref{prop.lapRV} and Corollary \ref{cor.lapRV2}). When
inserted in the spectral formula, the first term $R_\D$ at the right
hand side reproduces the solution to the free Dirac equation, and the
main part of our proof will be the estimate of second term
\begin{equation}\label{eq.Q}
  Q:=R_\D
  VR_\D(I+VR_\D)^{-1}.
\end{equation}
To this end, we shall need an explicit representation for
$R_{\D}(\lambda\pm i0)$, which is easily obtained from the formula
\begin{equation}\label{eq.RDR0}
  R_{\D}(\lambda\pm i0)=R_0(\lambda^2\pm i0)(\D+\lambda I_4).
\end{equation}
Recalling \eqref{eq.expl}, after an integration by parts we obtain
\begin{eqnarray}\label{eq.RD}
  R_{\D}(\lambda\pm i0)f & = & \frac{i\lambda}{4\pi}
  \int_{\R^3}\frac{e^{\pm i\lambda|x-y|}}{|x-y|}\left(
    I_4 \mp\sum_{j=1}^{3}\alpha_j\frac{x_j-y_j}{|x-y|}\right)f(y) d
  y\nonumber
  \\ \  & \  &
  +\frac{1}{4\pi} \int_{\R^3}
  \frac{e^{\pm i\lambda|x-y|}}{|x-y|^2}\sum_{j=1}^{3}
  \alpha_j\frac{x_j-y_j}{|x-y|}f(y) d y.
\end{eqnarray}
From here we derive immediately an analogous representation for
\begin{equation*}
  R_{\D}^2(\lambda)=
  \frac{\partial}{\partial\lambda}R_{\D}(\lambda);
\end{equation*}
indeed, differentiating \eqref{eq.RD} with respect to $\lambda$, we
get
\begin{eqnarray}\label{eq.RD2}
  R_{\D}^2(\lambda\pm i0)f & = & \frac{\lambda}{4\pi}
  \int_{\R^3}e^{\pm i\lambda|x-y|}
  \left(\mp
    I_4+\sum_{j=1}^{3}\alpha_j
    \frac{x_j-y_j}{|x-y|}\right)f(y)dy\nonumber \\ \  & \  &
  \pm\frac{i}{4\pi} \int_{\R^3}
  \frac{e^{\pm i\lambda|x-y|}}{|x-y|}\sum_{j=1}^{3}
  \alpha_j\frac{x_j-y_j}{|x-y|}f(y) d y.
\end{eqnarray}
We collect all the necessary estimates in the following lemma (we
write for simplicity $R_{\D}(\lambda)$ instead of $R_{\D}(\lambda\pm
i0)$ since the estimates are the same):

\begin{lemma}\label{lem.reD}
  Suppose that
  \begin{equation}\label{eq.assVbis}
    |V(x)|\leq\frac{C_{0}}{|x|^{1/2}\xx^{s}
      (|\log|x\|+1)^{\beta/2}},
  \end{equation}
  for some $s>\frac32$, $\beta>1$, $C_{0}>0$. Then the following
  estimates hold for all $\epsilon>0$ small enough and all
  $\lambda\in\mathbb{R}$:
  \begin{equation}\label{eq.VR2}
    \|\xx^{1/2+\epsilon}
    VR_{\D}^2(\lambda)f\|_{L^2}
    \leq C_{\epsilon}\cdot\xla\cdot
    \|\xx^{3/2+\epsilon} f\|_{L^2},
  \end{equation}
  \begin{equation}\label{eq.RVR}
    \|R_{\D}(\lambda)VR_{\D}(\lambda)f\|_{L^\infty}
    \leq
    C_{\epsilon}\cdot\xla^2\cdot
    \|\xx^{1/2+\epsilon}f\|_{L^2}
  \end{equation}
  \begin{equation}\label{eq.RVR2}
    \|R_{\D}^2(\lambda)VR_{\D}(\lambda)f\|_{L^\infty}
    +  \|R_{\D}(\lambda)VR_{\D}^2(\lambda)f\|_{L^\infty}
    \leq
    C_{\epsilon}\cdot\xla^2\cdot
    \|\xx^{3/2+\epsilon}f\|_{L^2}
  \end{equation}
  for some $C=C_\epsilon$ independent of $\lambda$.
\end{lemma}

\begin{proof}
  In the following we shall use the shorthand notation, for
  $s\in\mathbb{R}$,
  \begin{equation}\label{eq.L2g}
    \|f\|_{L^{2}_{\gamma}}:=\|\xx^{\gamma}f\|_{L^{2}}
  \end{equation}
  From the explicit representations \eqref{eq.RD} and \eqref{eq.RD2}
  we have the simple pointwise estimates
  \begin{equation}\label{eq.RDp}
    |R_{\D}(\lambda)f|\leq
    C(|\lambda|\cdot|x|^{-1}+|x|^{-2})* f,\qquad
    |R_{\D}^2(\lambda)f|\leq
    C(|\lambda|+|x|^{-1})* f.
  \end{equation}
  Since $|x|^{-1}\in L^{3,\infty}$, by the Young inequality in Lorentz
  spaces (see the Appendix) we get
  \begin{align*}
    \|VR_{\D}^2(\lambda)f\|_{L^2_\gamma} & \leq
    \|V\|_{L^2_\gamma}\cdot|\lambda|\cdot \|1*f\|_{L^\infty}+
    \|V\|_{L^2_\gamma}\||x|^{-1}*f\|_{L^\infty} \\ \ & \leq
    \|V\|_{L^2_\gamma} \left(
      |\lambda|\cdot\|f\|_{L^1}+\|f\|_{L^{3/2,1}}\right).
  \end{align*}
  By the obvious inequalities valid for all $\epsilon>0$
  \begin{equation}\label{eq.lor}
    \|f\|_{L^{1}}\leq C(\epsilon)
    \|f\|_{L^{2}_{3/2+\epsilon}},\qquad
    \|f\|_{L^{3/2,1}}\leq C(\epsilon)
    \|f\|_{L^{2}_{1/2+\epsilon}},\qquad
  \end{equation}
  we arrive at the first estimate
  \begin{equation}\label{eq.int1}
    \|VR_{\D}^2(\lambda)f\|_{L^2_\gamma}\leq
    C(\epsilon) \|V\|_{L^2_\gamma}\xla
    \|f\|_{L^{2}_{3/2+\epsilon}}.
  \end{equation}
  Since $\|V\|_{L^2_\gamma}<\infty$ by assumption \eqref{eq.assVbis}
  as soon as $\gamma=1/2+\epsilon<s-1$, we see that \eqref{eq.VR2}
  follows provided $\epsilon$ is suitably small.

  In a similar way, in order to prove \eqref{eq.RVR} we use again
  \eqref{eq.RDp} and we write (recall that $|x|^{-2}\in
  L^{3/2,\infty}$)
  \begin{align*}
    \|R_{\D}(\lambda)VR_{\D}(\lambda)f\|_{L^\infty} & \leq C \left(
      |\lambda|\cdot\||x|^{-1}*VR_{\D}f\|_{L^{\infty}}+
      \||x|^{-2}*VR_{\D}f\|_{L^{\infty}}
    \right)    \\
    &\leq C\left( |\lambda|\cdot\|VR_{\D}(\lambda)f\|_{L^{3/2,1}}+
      \|VR_{\D}(\lambda)f\|_{L^{3,1}} \right) .
  \end{align*}
  For the first term we can write, recalling again \eqref{eq.RDp},
  \begin{align}\label{eq.frst}
    \|VR_{\D}(\lambda)f\|_{L^{3/2,1}}&\leq \|V\|_{L^{3/2,1}}
    |\lambda|\cdot \||x|^{-1}*f\|_{L^{\infty}}
    +\|V\|_{L^{2}}\||x|^{-2}*f\|_{L^{6,2}}
    \\
    &\leq \|V\|_{L^{3/2,1}} |\lambda|\cdot \|f\|_{L^{3/2,1}}
    +\|V\|_{L^{2}}\|f\|_{L^{2}}\nonumber
    \\
    &\leq\left( \|V\|_{L^{3/2,1}} |\lambda| +\|V\|_{L^{2}}\right)
    \|f\|_{L^{2}_{3/2+\epsilon}}\nonumber
  \end{align}
  (see \eqref{eq.lor}), while for the second term we have
  \begin{align}\label{eq.snd}
    \|VR_{\D}(\lambda)f\|_{L^{3,1}}&\leq \|V\|_{L^{3,1}}
    |\lambda|\cdot \||x|^{-1}*f\|_{L^{\infty}}
    +\|V\|_{L^{6,2}}\||x|^{-2}*f\|_{L^{6,2}}
    \\
    &\leq \|V\|_{L^{3,1}} |\lambda|\cdot \|f\|_{L^{3/2,1}}
    +\|V\|_{L^{6,2}}\|f\|_{L^{2}}\nonumber
    \\
    &\leq\left( \|V\|_{L^{3,1}} |\lambda| +\|V\|_{L^{6,2}}\right)
    \|f\|_{L^{2}_{3/2+\epsilon}}\nonumber
  \end{align}
  where we have used \eqref{eq.lor} and the trivial inequality
  $\|f\|_{L^{2}}\leq\|f\|_{L^{2}_{\gamma}}$, $\forall\gamma>0$.
  Summing up, we get
  \begin{equation}\label{eq.thrd}
    \|R_{\D}(\lambda)VR_{\D}(\lambda)f\|_{L^\infty}  \leq
    C\cdot C(V)\xla^2 \|f\|_{L^{2}_{3/2+\epsilon}}
  \end{equation}
  where the quantity
  \begin{equation}\label{eq.hypV}
    C(V):=\|V\|_{L^{3/2,1}}+
    \|V\|_{L^{3,1}}+
    \|V\|_{L^{6,2}}+
    \|V\|_{L^{2}}<\infty
  \end{equation}
  is finite by assumption \eqref{eq.assVbis} (see also the Appendix \ref{appendix}).

  The proof of \eqref{eq.RVR2} is similar: by \eqref{eq.RDp} we get
  \begin{align*}
    \|R_{\D}^{2}(\lambda)VR_{\D}(\lambda)f\|_{L^\infty} & \leq C
    \left( |\lambda|\cdot\|1*VR_{\D}f\|_{L^{\infty}}+
      \||x|^{-1}*VR_{\D}f\|_{L^{\infty}}
    \right)    \\
    &\leq C\left( |\lambda|\cdot\|VR_{\D}(\lambda)f\|_{L^{1}}+
      \|VR_{\D}(\lambda)f\|_{L^{3/2,1}} \right) .
  \end{align*}
  We have already estimated the second term in \eqref{eq.frst}, and
  for the first one we have
  \begin{align}\label{eq.frt}
    \|VR_{\D}(\lambda)f\|_{L^{1}}&\leq \|V\|_{L^{3/2}} |\lambda|\cdot
    \||x|^{-1}*f\|_{L^{3}} +\|V\|_{L^{3}}\||x|^{-2}*f\|_{L^{3/2}}
    \\
    &\leq\left( \|V\|_{L^{3/2}} |\lambda| +\|V\|_{L^{3}}\right)
    \|f\|_{L^{1}}\nonumber
    \\
    &\leq\left( \|V\|_{L^{3/2}} |\lambda| +\|V\|_{L^{3}}\right)
    \|f\|_{L^{2}_{3/2+\epsilon}}\nonumber
  \end{align}
  and hence
  \begin{equation}\label{eq.thrd2}
    \|R_{\D}^{2}(\lambda)VR_{\D}(\lambda)f\|_{L^\infty}  \leq
    C\cdot C'(V)\xla^2 \|f\|_{L^{2}_{3/2+\epsilon}}
  \end{equation}
  where the quantity
  \begin{equation}\label{eq.hypV2}
    C'(V):=\|V\|_{L^{3/2}}+
    \|V\|_{L^{3/2,1}}+
    \|V\|_{L^{3}}+
    \|V\|_{L^{2}}<\infty
  \end{equation}
  is finite again by assumption \eqref{eq.assVbis}.

  Finally, the last estimate can be obtained as follows:
  \begin{align*}
    \|R_{\D}(\lambda)VR_{\D}^{2}(\lambda)f\|_{L^\infty} & \leq C
    \left( |\lambda|\cdot\||x|^{-1}*VR_{\D}^{2}f\|_{L^{\infty}}+
      \||x|^{-2}*VR_{\D}^{2}f\|_{L^{\infty}}
    \right)    \\
    &\leq C\left( |\lambda|\cdot\|VR_{\D}(\lambda)f\|_{L^{3/2,1}}+
      \|VR_{\D}(\lambda)f\|_{L^{3,1}} \right) .
  \end{align*}
  Proceeding as above, we estimate
  \begin{align}\label{eq.fif}
    \|VR_{\D}^{2}(\lambda)f\|_{L^{3/2,1}}&\leq \|V\|_{L^{3/2,1}}
    |\lambda|\cdot \|1*f\|_{L^{\infty}}
    +\|V\|_{L^{3/2,1}}\||x|^{-1}*f\|_{L^{\infty}}
    \\
    &\leq \|V\|_{L^{3/2,1}}4\xla \left(\|f\|_{L^{1}}+
      \|f\|_{L^{3/2,1}}\right) \nonumber
    \\
    &\leq \|V\|_{L^{3/2,1}}4\xla \|f\|_{L^2_{3/2+\epsilon}} \nonumber
  \end{align}
  and
  \begin{align}\label{eq.sxt}
    \|VR_{\D}^{2}(\lambda)f\|_{L^{3,1}}&\leq \|V\|_{L^{3,1}}
    |\lambda|\cdot \|1*f\|_{L^{\infty}}
    +\|V\|_{L^{3,1}}\||x|^{-1}*f\|_{L^{\infty}}
    \\
    &\leq \|V\|_{L^{3,1}}4\xla \left(\|f\|_{L^{1}}+
      \|f\|_{L^{3/2,1}}\right)\nonumber
    \\
    &\leq \|V\|_{L^{3,1}}4\xla \|f\|_{L^2_{3/2+\epsilon}}\nonumber
  \end{align}
  whence
  \begin{equation}\label{eq.sept}
    \|R_{\D}(\lambda)VR_{\D}^{2}(\lambda)f\|_{L^\infty}  \leq
    C\cdot C''(V)\xla^2 \|f\|_{L^{2}_{3/2+\epsilon}}
  \end{equation}
  where the quantity
  \begin{equation}\label{eq.hypV3}
    C''(V):=
    \|V\|_{L^{3/2,1}}+
    \|V\|_{L^{3,1}}<\infty
  \end{equation}
  is finite by assumption \eqref{eq.assVbis}.
\end{proof}

\begin{remark}\label{rem.spec2}
  The same remark concerning the simpler version of the spectral
  formula \eqref{eq.spec1} and the integration by parts formula
  \eqref{eq.intp} applies also to the Dirac resolvent, with obvious
  modifications in the proof.
\end{remark}

\section{Proof of Theorem \ref{th.1}}\label{sec.onde}

Let $\left(\varphi_j\right)_{j=0,1,\dots}$ be a standard
Paley-Littlewood partition of the unity, with the properties
\begin{equation}\label{eq.PL}
  \varphi_j(\lambda)=\varphi_0(2^{-j}\lambda), \qquad
  \varphi_0+\sum_{j\geq1}\varphi_j=1,
\end{equation}
for a suitable $\varphi_0\in\mathcal C^\infty_0$.  We consider the
Cauchy problem
\begin{equation}\label{eq.ondej}
  \left\{\begin{array}{l}u_{tt}(t,x)-\Delta u(t,x)+W(x,D)u=0 \\
      u(0,x)=0,\quad
      u_t(0,x)=\varphi_j(\sqrt{-\Delta+W})g(x),
    \end{array}\right.
\end{equation}
The solution can be represented using the spectral formula as follows:
\begin{equation}\label{eq.solonde}
  u(t,x)=\frac{1}{2\pi i}
  \int_{0}^{+\infty}
  \varphi_j(\sqrt\lambda)
  \frac {\sin(t\sqrt\lambda)}{\sqrt\lambda}
  R(\lambda) gd\lambda,
\end{equation}
and after an integration by parts (see Remark \ref{rem.spec1}) this
gives
\begin{equation}\label{eq.solonde2}
  u(t,x)=\frac C t
  \int_{0}^{+\infty}
  \cos(t\sqrt\lambda)
  \left[\partial_{\lambda}
    \varphi_j(\sqrt\lambda) R(\lambda) g
    +
    \varphi_j(\sqrt\lambda)
    \partial_{\lambda}R(\lambda) g
  \right]d\lambda.
\end{equation}
Thus, recalling estimates \eqref{eq.stimaRW} and \eqref{eq.stimaRW2},
we have
\begin{equation*}
  |u(t,x)|\leq\frac{C}{t}
  \|\xx w_{\beta}^{1/2}g\|_{L^2}
  \int_{0}^{+\infty}
  \left(
    |\partial_{\lambda}\varphi_j(\sqrt\lambda)|
    +\left(1+\frac{1}{\sqrt\lambda}\right)
    |\varphi_j(\sqrt\lambda)|
  \right) d\lambda
\end{equation*}
and a change of variables $\lambda=2^{2j}\mu$ in the integral gives
immediately
\begin{equation}\label{eq.estj}
  |u(t,x)|\leq \frac{C}{t}2^{2j}
  \|\xx w_{\beta}^{1/2}g\|_{L^2}
\end{equation}
with some constant $C$ independent of $j$ and $g$.

If we now define as usual
\begin{equation*}
    \widetilde\varphi_{j}=
    \varphi_{j-1}+\varphi_{j}+\varphi_{j+1},\quad
    \varphi_{-1}=0,
\end{equation*}
so that $\varphi_{j}\equiv\varphi_{j}\widetilde\varphi_{j}$,
we see that the Cauchy problem \eqref{eq.ondej} can be
written equivalently
\begin{equation}\label{eq.ondejtilde}
  \left\{\begin{array}{l}u_{tt}(t,x)-\Delta u(t,x)+W(x,D)u=0 \\
      u(0,x)=0,\quad
      u_t(0,x)=\varphi_j(\sqrt{-\Delta_W})
      \widetilde\varphi_{j}(\sqrt{-\Delta_W})g(x),
    \end{array}\right.
\end{equation}
hence our estimate \eqref{eq.estj} implies also the estimate
\begin{equation}\label{eq.estjj}
  |u(t,x)|\leq \frac{C}{t}2^{2j}
  \|\xx w_{\beta}^{1/2}
  \widetilde\varphi_{j}(\sqrt{-\Delta_W})g\|_{L^2}.
\end{equation}
Finally, consider the original Cauchy problem \eqref{eq.onde}, and
decompose $g$ as a sum
$g=\sum_{j\geq0}\varphi_{j}(\sqrt{-\Delta_W})g(x)$.
By estimate \eqref{eq.estjj} we obtain easily estimate
\eqref{eq.decwave}.
\begin{equation}\label{eq.estlast}
  |u(t,x)|\leq \frac{C}{t}\sum_{j\geq0}2^{2j}
  \|\xx w_{\beta}^{1/2}\varphi_{j}(\sqrt{-\Delta_W})g\|
  _{L^{2}}.
\end{equation}
The computations in the case of initial data of the form
\begin{equation*}
    u(0,x)=f,\qquad u_{t}(0,x)=0
\end{equation*}
are completely analogous, and we thus obtain estimate
\eqref{eq.decwavef}.

\begin{remark}\label{rem.foll}
In view of the application to the Dirac system, the following remark
will be useful. If the initial datum $g$ has the form
\begin{equation}\label{eq.Ds}
    g=(-\Delta_{W})^{s}h
\end{equation}
for some $s>0$, a direct application of estimate \eqref{eq.estlast}
would give only
\begin{equation}\label{eq.estlasts}
  |u(t,x)|\leq \frac{C}{t}\sum_{j\geq0}2^{2j}
  \|\xx w_{\beta}^{1/2}\varphi_{j}(\sqrt{-\Delta_W})
 ( -\Delta_{W})^{s}h\|
  _{L^{2}}.
\end{equation}
Actually, if we go back to the spectral formula \eqref{eq.solonde2},
we see that the solution can be written
\begin{equation}\label{eq.solonde2s}
  u(t,x)=\frac C t
  \int_{0}^{+\infty}\lambda^{s/2}
  \cos(t\sqrt\lambda)
  \left[\partial_{\lambda}
    \varphi_j(\sqrt\lambda) R(\lambda) h
    +
    \varphi_j(\sqrt\lambda)
    \partial_{\lambda}R(\lambda) h
  \right]d\lambda.
\end{equation}
with an additional factor $\lambda^{s/2}$. Thus, proceeding as
above, we arrive at the simpler estimate
\begin{equation}\label{eq.estlastgarr}
  |u(t,x)|\leq \frac{C}{t}\sum_{j\geq0}2^{(2+s)j}
  \|\xx w_{\beta}^{1/2}\varphi_{j}(\sqrt{-\Delta_W})h\|
  _{L^{2}}.
\end{equation}

\end{remark}

We now prove estimate \eqref{eq.decwaveH2}
under the stronger assumption \eqref{eq.ipfinregW}
on the potential $W(x,D)$. Consider first the case of initial data
of the form
\begin{equation*}
        u(0,x)=0,\qquad u_{t}(0,x)=g.
\end{equation*}
We can write $g$ as follows:
\begin{equation*}
    g=(1-\Delta+W)^{-1-\epsilon}(1-\Delta+W)^{1+\epsilon}g
\end{equation*}
for some fixed $\epsilon>0$. Then the solution $u$ can be
represented as
\begin{equation*}
      u(t,x)=\frac{1}{2\pi i}
  \int_{0}^{+\infty}
  \psi(\sqrt\lambda)
  \frac {\sin(t\sqrt\lambda)}{\sqrt\lambda}
  R(\lambda) hd\lambda
\end{equation*}
where
\begin{equation*}
    h=(1-\Delta+W)^{1+\epsilon}g,\qquad
    \psi(\sqrt{\lambda})=(1+\lambda)^{1+\epsilon}.
\end{equation*}
Proceeding as above, after an integration by parts we arrive at
\begin{equation*}
    |u(t,x)|\leq
    \frac C t\|\xx w_{\beta}^{1/2}h\|_{L^{2}}
    \int_{0}^{+\infty}((1+\lambda)^{-1-\epsilon}+
        (1+\lambda)^{-2-\epsilon})d\lambda
\end{equation*}
and hence
\begin{equation}\label{eq.estdelta}
    |u(t,x)|\leq
    \frac C t\|\xx w_{\beta}^{1/2}
        (1-\Delta+W)^{1+\epsilon}g\|_{L^{2}}\leq
        \frac C t\|\xx^{3/2+\epsilon}
        (1-\Delta+W)^{1+\epsilon}g\|_{L^{2}}.
\end{equation}
To conclude the proof of the Theorem, it remains to show that
\begin{equation}\label{eq.sobin}
    \|\xx^{3/2+\epsilon}
           (1-\Delta+W)^{1+\epsilon}g\|_{L^{2}}\leq
    \|\xx^{3/2+\epsilon}g\|_{H^{2+2\epsilon}}.
\end{equation}
We start from the inequality
\begin{equation*}
    \|\xx^{s}(1-\Delta+W)f\|_{L^{2}}\leq
     \|\xx^{s}f\|_{H^{2}}
\end{equation*}
which is obviously valid for any $s\geq 0$. By a standard
complex interpolation argument, interpolating with the
trivial inequality
\begin{equation*}
    \|\xx^{s}f\|_{L^{2}}\leq
     \|\xx^{s}f\|_{L^{2}}
\end{equation*}
we obtain that
\begin{equation*}
    \|\xx^{s}(1-\Delta+W)^{\epsilon}f\|_{L^{2}}\leq
     \|\xx^{s}f\|_{H^{2\epsilon}}
\end{equation*}
for all $0\leq \epsilon\leq1$ and all $s\geq0$.
This implies
\begin{equation}\label{eq.quasiult}
    \|\xx^{s}(1-\Delta+W)^{1+\epsilon}f\|_{L^{2}}\leq
     \|\xx^{s}(1-\Delta+W)f\|_{H^{2\epsilon}}\leq
    \|\xx^{s}f\|_{H^{2+2\epsilon}}+
    \|\xx^{s}Wf\|_{H^{2\epsilon}}.
\end{equation}
The last term is of the form
\begin{equation}\label{eq.pon1}
    \|\xx^{s}W(x,D)f\|_{H^{2\epsilon}}\leq
     \|\xx^{s}a(x)Df\|_{H^{2\epsilon}}+
     \|\xx^{s}b(x)f\|_{H^{2\epsilon}};
\end{equation}
in order to estimate it, we recall the Kato-Ponce inequality
(see \cite{KP})
\begin{equation}\label{eq.KP}
    \|\DD^{q}(vw)\|_{L^{p}}\leq
    C\|\DD^{q}v\|_{L^{p_{1}}}\|w\|_{L^{p_{2}}}+
    C\|v\|_{L^{p_{3}}}\|\DD^{q}w\|_{L^{p_{4}}}
\end{equation}
which is valid for all $q\geq0$, $p^{-1}=p_{1}^{1}+p_2^{-1}=
p^{-1}_{3}+p^{-1}_{4}$. With the choices $v(x)=a(x)$,
$w(x)=\xx^{s}Df(x)$, $q=2\epsilon$, $p_{1}=p_{3}=\infty$
and $p_{2}=p_{4}=2$, we obtain
\begin{equation*}
    \|\DD^{2\epsilon}\xx^{s}a(x)Df\|_{L^{2}}\leq
    C\|\DD^{2\epsilon}a\|_{L^{\infty}}\|\xx^{s}Df\|_{L^{2}}
    +C
    \|a\|_{L^{\infty}}\|\DD^{2\epsilon}(\xx^{s}Df)\|_{L^{2}}.
\end{equation*}
Now it is clear that
\begin{equation*}
    \|\DD^{2\epsilon}(\xx^{s}Df)\|_{L^{2}}\leq
    C \|\xx^{s}f\|_{H^{1+2\epsilon}}
\end{equation*}
(use again complex interpolation between the cases $\epsilon=0$
and $\epsilon=1$) and in conclusion we obtain
\begin{equation*}
    \|\DD^{2\epsilon}\xx^{s}a(x)Df\|_{L^{2}}\leq
    C
    \|\DD^{2\epsilon}a\|_{L^{\infty}}
        \|\xx^{s}f\|_{H^{1+2\epsilon}}.
\end{equation*}
Here we have used the simple fact that
\begin{equation*}
     \|a\|_{L^{\infty}}\leq
     C \|\DD^{2\epsilon}a\|_{L^{\infty}}.
\end{equation*}
The corresponding estimate for the electric term is analogous
(actually simpler):
\begin{equation*}
    \|\DD^{2\epsilon}\xx^{s}b(x)f\|_{L^{2}}\leq
    C
    \|\DD^{2\epsilon}b\|_{L^{\infty}}
            \|\xx^{s}f\|_{H^{2\epsilon}}.
\end{equation*}
Recalling now \eqref{eq.quasiult} and \eqref{eq.pon1} we conclude
the proof of estimate \eqref{eq.decwaveH2}.

On the other hand, when the data are of the form
\begin{equation*}
    u(0,x)=f,\qquad u_{t}(0,x)=0
\end{equation*}
the computations are completely analogous and we obtain estimate
\eqref{eq.decwaveH2f}
under the stronger assumptions
\eqref{eq.ipfinregWf} on the coefficients.

\section{Proof of Theorem \ref{th.2}}\label{sec.th2}

\begin{remark}\label{rem.reduction}
We notice that Theorem \ref{th.1} (and Remark \ref{rem.fg})
can be trivially
extended to a \emph{system} of wave equations of the form
\begin{equation}\label{eq.systw}
    u_{tt}-(\nabla+iA(x))^{2}u+
        B(x)u=0
\end{equation}
where $u(t,x)$ is a $\mathbb{C}^{N}$ valued function
and $A_{1}(x)$, $A_{2}(x)$, $A_{3}(x)$, $B(x)$ are $\mathbb{C}^{N\times N}$
matrices whose coefficients satisfy the assumptions of the
Theorem. The resulting dispersive estimates have
exactly the same form as in the scalar case.
\end{remark}

Consider now the Cauchy problem
\begin{equation}\label{dirac.pot2}
  \left\{\begin{array}{l}iu_t-\D u-V(x)u=0 \\
      u(0,x)=f(x).\end{array}\right.
\end{equation}
If we apply to the pertubed Dirac system the operator
$i\partial_{t}+\D+V$ we obtain that $u$ is also a solution of a
4$\times$4 system of perturbed wave equations of the form
\eqref{eq.systw} with
\begin{equation}\label{eq.coeff}
    A_{j}(x)=-\frac12(\alpha_{j}V(x)+V(x)\alpha_{j}),
\end{equation}
\begin{equation}\label{eq.coeff2}
    B(x)=\D V(x)+V(x)^{2}+A_{1}^{2}+A_{2}^{2}+A_{3}^{2}
    +i\sum\partial_{j}A_{j}
\end{equation}
and initial data
\begin{equation}\label{eq.datau}
    u(0,x)=f,\qquad
    u_{t}(0,x)=i^{-1}(\D+V)f.
\end{equation}
Note that the perturbed operator
\begin{equation}\label{eq.DW}
    -\Delta_{W}=-(\nabla+iA(x))^{2}+
        B(x)
\end{equation}
is exactly the square of the operator $\D+V$:
\begin{equation}\label{eq.squ}
     -\Delta_{W}=(\D+V)^{2}
\end{equation}
and hence the initial data for \eqref{eq.systw} can be written
\begin{equation}\label{eq.datau2}
    u(0,x)=f,\qquad
    u_{t}(0,x)=i^{-1}(-\Delta_{W})^{1/2}f.
\end{equation}
We are in position to apply to the solution $u$ the estimates
already proved in Theorem \ref{th.1}; keeping Remark \ref{rem.foll}
into account, we arrive easily at the estimate
\begin{equation}\label{eq.finD}
  |u(t,x)|\leq \frac{C}{t}\sum_{j\geq0}2^{3j}
  \|\xx w_{\beta}^{1/2}\varphi_{j}(\D+V)f\|
  _{L^{2}},
\end{equation}
provided the coefficients $a_{j}(x)$ and $b(x)$ satisfy the assumptions
\eqref{eq.hpWfin}. Recalling the explicit form \eqref{eq.coeff}
of the coefficients in terms of $V(x)$, we see that $V$ must satisfy
the conditions
\begin{equation*}
    |V(x)|\leq \frac{C_{0}}{|x|\xx(|\log|x||+1)^{\beta}}
\end{equation*}
from the magnetic term, and
\begin{equation*}
        |V(x)^{2}|+|DV(x)|\leq \frac{C_{0}}
        {|x|^{2}(|\log|x||+1)^{\beta}},\qquad
\end{equation*}
from the electric term, for some $\beta>1$ and some small constant
$C_{0}$. Summing up, we obtain that \eqref{eq.finD} holds under
assumption \eqref{ipfinVreg}.

The estimate in terms of the Sobolev norm can be
obtained in exactly the same way as for the
perturbed wave equation. Indeed, proceeding as in \eqref{eq.estdelta}
we arrive at the estimate
\begin{equation}\label{eq.finD2}
    |u(t,x)|\leq\frac Ct
    \|\xx^{3/2+\epsilon}
    (-\Delta_{W})^{3/2+\epsilon}f\|_{L^{2}}.
\end{equation}
The same arguments used at the end of Section \ref{sec.onde}
give here
\begin{equation}\label{eq.finD3}
    |u(t,x)|\leq\frac Ct
    \|\xx^{3/2+\epsilon}
           f\|_{H^{3+2\epsilon}}
\end{equation}
provided
\begin{equation*}
    \DD^{1+2\epsilon}A_{j}\in L^{\infty},\quad
    \DD^{1+2\epsilon}B\in L^{\infty},
\end{equation*}
which is implied by
\begin{equation*}
    \DD^{2+2\epsilon}V\in L^{\infty}.
\end{equation*}

\section{Proof of Theorem \ref{th.3}}\label{sec.decdirac}

By exploiting the connection between the massless Dirac and the
wave equation, it is easy to obtain an optimal dispersive
estimate in the unperturbed case. Indeed, let $u(t,x)$ be
a smooth solution of the free massless Dirac equation
\begin{equation}\label{eq.dirac2}
  iu_t(t,x)=\D u(t,x)
\end{equation}
with initial data
\begin{equation}\label{eq.diracd}
    u(0,x)=f(x).
\end{equation}
Recall now the identity
\begin{equation*}
     (i\partial_{t}+\D)(i\partial_{t}-\D)
     =(\Delta-\partial^{2}_{tt}) I_{4};
\end{equation*}
if we apply the operator $i\partial_{t}+\D$ to the system
\eqref{eq.dirac2} we immediately obtain that $u$ solves the
Cauchy problem for the wave equation
\begin{equation*}
    u_{tt}-\Delta u=0
\end{equation*}
with initial data
\begin{equation*}
    u(0,x)=f,\qquad
    u_{t}(0,x)=i^{-1} \D f.
\end{equation*}
Then, as a consequence of the well known decay estimates for
solutions to the free wave equation (see e.g.
\cite{SS}), we obtain
\begin{equation*}
    |u(t,x)|\leq \frac C t \left(\|f\|_{\dot B^{2}_{1,1}}
      +\|Df\|_{\dot B^{1}_{1,1}}\right)
\end{equation*}
and hence
\begin{equation}\label{eq.decayu}
        |u(t,x)|\leq \frac C t \|f\|_{\dot B^{2}_{1,1}}.
\end{equation}
Here $\dot B^{s}_{1,1}$ is the homogeneous Besov space,
with norm
\begin{equation*}
    \|v\|_{\dot B^{s}_{1,1}}=\sum_{j\in\mathbb{Z}}
      2^{js}
        \|\phi_{j}(\sqrt{-\Delta})v\|_{L^{1}}
\end{equation*}
where $\phi_{j}$ now is a \emph{homogeneous} Paley-Littlewood sequence,
i.e., fixed a test function $\psi(r)\in C^{\infty}_{0}$ such that $\psi(r)=1$ for $r<1$, $\psi(r)=0$ for $r>2$, we have
$\phi_{j}(r)=\psi(2^{-j+2}r)-\psi(2^{-j+1}r)$ for all
$j\in\mathbb{Z}$.

The proof of Theorem \ref{th.3} follows the same lines as the proof of
Theorem \ref{th.1}. Consider the Cauchy problem
with frequency truncated data
\begin{equation}\label{eq.diracomp}
  \left\{\begin{array}{l}
  iu_t(t,x)=\D_V u(t,x) \\ u(0,x)=\varphi_j
  (\D_V)f,\end{array}\right.
\end{equation}
where $(\varphi_j(\lambda))_{j=0,1,\dots}$ is the standard
Paley-Littlewood partition of the unity defined in
\eqref{eq.PL}. By means of spectral formula, we can represent
the solution of \eqref{eq.diracomp} as
\begin{equation}\label{eq.dirspec}
    u(t,x)=\frac 1{2\pi i}\int_{-\infty}^{+\infty}\varphi_j(\lambda)
    e^{i\lambda t}\mathcal I[R_V(\lambda)]f\,d\lambda.
\end{equation}
Using the identity
\begin{equation}\label{eq.resdirac}
  R_V(\lambda)=R_{\D}-R_{\D}VR_{\D}(I+VR_{\D})^{-1},
\end{equation}
which is valid thanks to Corollary \ref{cor.lapRV2},
we can split the integrals in \eqref{eq.dirspec} into two
terms, the first one containing the contribution of the free
resolvent $R_{\D}$ and the second one containing the contribution of
the operator $R_{\D}VR_{\D}(I+VR_{\D})^{-1}$. The
first term
\begin{equation*}
  A:=\frac{1}{2\pi
    i}\int_{-\infty}^{+\infty}\varphi_j(\lambda)e^{i\lambda
    t}\Im \left[R_{\D}(\lambda)\right]f\diff\lambda
\end{equation*}
was estimated above (see \eqref{eq.decayu}); it remains to
estimate the term
\begin{equation}\label{eq.pezzoB}
  B=-\frac{1}{2\pi
    i}\int_{-\infty}^{+\infty}\varphi_j(\lambda)e^{i\lambda
    t}\Im\left[Q(\lambda)\right]f\diff\lambda,
\end{equation}
where
\begin{equation*}
  Q(\lambda):=R_{\D}(\lambda)VR_{\D}(\lambda)
  (I+VR_{\D}(\lambda))^{-1}.
\end{equation*}
After an integration by parts, we obtain
\begin{equation}\label{eq.dirpart}
  B=-\frac{1}{2\pi t}\left[\int_{-\infty}^{+\infty}\varphi_j(\lambda)
    e^{i\lambda t}\frac{\partial}{\partial\lambda}
    \mathcal{I}(Q(\lambda))f\diff\lambda+\int_{-\infty}^{+\infty}
    \varphi'_j(\lambda)e^{i\lambda
      t}\mathcal{I}[Q(\lambda)]f\diff\lambda\right];
\end{equation}
an explicit computation shows that
\begin{eqnarray*}
  \frac{\partial Q}{\partial\lambda} & = &
  R_{\D}^2VR_{\D}(I_4+VR_{\D})^{-1}+
  R_{\D}VR_{\D}^2(I_4+VR_{\D})^{-1} \\ \
  & \  &
  +R_{\D}VR_{\D}(I_4+VR_{\D})^{-1}VR_{\D}^2
  (I_4+VR_{\D})^{-1}.
\end{eqnarray*}
Now we can apply
 Lemma \ref{lem.reD}: under assumption \eqref{eq.ipfinV},
estimates \eqref{eq.VR2}, \eqref{eq.RVR} and \eqref{eq.RVR2}
are satisfied, and the Lemma gives
\begin{equation}\label{eq.stimaQ}
  \|Q(\lambda)f\|_{L^\infty} \leq
  C\xla^2\|\xx^{3/2+\epsilon}f\|_{L^2},
\end{equation}
\begin{equation}\label{eq.stimaQ'}
  \left|\left|\frac{\partial}{\partial\lambda}
    Q(\lambda)f\right|\right|_{L^\infty} \leq
  C\xla^3\|\xx^{3/2+\epsilon}f\|_{L^2},
\end{equation}
for some $C>0$. Using \eqref{eq.stimaQ} and \eqref{eq.stimaQ'} in
\eqref{eq.dirpart} we arrive at the estimate
\begin{equation*}
  |B|\leq\frac
    Ct\|f\|_{L^2_{3/2+\epsilon}}\left[\int_{-\infty}^{+\infty}\left(
      \xla^3|\varphi_j(\lambda)|+\xla^2\left|
        \varphi'_j(\lambda)
      \right|\right)\diff\lambda\right].
\end{equation*}
Recalling that $\phi_{j}(\lambda)=\phi_{0}(2^{-j}\lambda)$,
after a change of variables $2^{-j}\lambda=\mu$ we easily
obtain
\begin{equation}\label{stimaB}
  |B|\leq\frac Ct2^{4j}
  \|\xx^{3/2+\epsilon}f\|_{L^2}.
\end{equation}
From this point on, we can proceed as in the proof of Theorem
\ref{th.1} and complete the proof of Theorem \ref{th.3}.

\appendix
\section{Lorentz spaces}\label{appendix}

For the convenience of the reader, we recall here
the definitions and the main
properties of the Lorentz spaces $L^{p,q}$, in view
of the applications needed in the proof of our results.

For any  measurable function $f:\mathbb R^n\to\mathbb C$
and any
$s\geq0$ we define the \emph{upper-level} $E^f_s$ as the set
$$E^f_s:=\{x:|f(x)|>s\}.$$
The \emph{non-increasing rearrangement} of $f$ is
then the function
$$f^*(t):=\inf\{s>0:|E^f_s|\leq t\},\qquad
    t\in(0,+\infty).$$
It is also useful to consider the average of $f^*$ defined by
$$f^{**}(t)=\frac 1t\int_0^tf^*(r)\,dr.$$

The standard definition of the Lorentz spaces is
the following:
\begin{definition}\label{def.lorentz}
    For any $1\leq p<\infty$ and $1\leq q\leq\infty$
    we define the quasinorm $\|f\|_{L^{p,q}}$ as
    follows:
    \begin{equation}\label{norm.lorentz}
    \|f\|_{L^{p,q}}=
    \begin{cases}
    \left[\int_0^\infty(t^{1/p}f^{*}(t))^q\frac{dt}{t}
    \right]^{1/q},  &  1\leq q<\infty
    \\
    \sup_{t>0}t^{1/p}f^{*}(t),
    &   q=\infty.
    \end{cases}
    \end{equation}
    When $p\neq1$, if we replace $f^{*}$ with $f^{**}$ in the
    above definitions we obtain an equivalent quasinorm which
    is actually a norm (see \cite{BL}, \cite{cald}).
    The \emph{Lorentz space} $L^{p,q}$ is defined by
    \begin{equation}\label{lorentz.space}
        L^{p,q}=\{f:\|f\|_{L^{p,q}}<\infty\}.
    \end{equation}
    Moreover we define
    $$L^{1,1}:=L^{1},\qquad L^{\infty,\infty}=L^{\infty}.$$
    The spaces $L^{\infty,q}$ for $1\leq q<\infty$ are usually
    left undefined (although $L^{\infty,1}$ is defined in \cite{cald}
    as the closure of $L^{\infty}$ compactly supported functions in the
    $L^{\infty}$ norm).
\end{definition}

With the above definitions, one obtains the elementary
properties
\begin{equation*}
    L^{p,p}=L^{p},\qquad 1\leq p\leq\infty;
\end{equation*}
\begin{equation*}
    L^{p,q_{1}}\subseteq L^{p,q_{2}},\qquad
    1<p<\infty,\quad
    1\leq q_{1}\leq q_{2}\leq\infty
\end{equation*}
(with continuous embedding). When the second index is $\infty$
we obtain the weak Lebesge spaces (Marcinkiewicz spaces):
\begin{equation*}
    L^{p,\infty}=L^{p}_{w},\qquad 1\leq p\leq\infty.
\end{equation*}
Moreover, the Lorentz spaces can
be obtained by an equivalent construction using real interpolation:
\begin{equation*}
    L^{p,q}=(L^{p_{0}},L^{p_{1}})_{\theta,q},\qquad
    p^{-1}=(1-\theta)p_{0}^{-1}+\theta p_{1}^{-1}
\end{equation*}
provided
\begin{equation*}
   p_{0}<p_{1},\qquad p_{0}<q\leq\infty,\qquad 0<\theta<1.
\end{equation*}

An alternative characterization of the Lorentz norm can be given
using the so-called \emph{atomic decomposition}:

\begin{lemma}\label{lorentz.atomic}
    Let $f:\mathbb R^n\to\mathbb C$ be a measurable function
    and let $1\leq p<\infty,\ 1\leq q\leq\infty$;
    then $f\in L^{p,q}$ if and only if there exist a sequnce
    of sets $(E_j)_{j\in\mathbb Z}$ and a sequence of numbers
    $a=(a_j)_{j\in\mathbb Z}$ such that $|E_j|=O(2^j)$, $a\in l^q$
    and the following estimate
    \begin{equation}\label{eq.lorentz.atomic}
        |f(x)|\leq C\sum_{j\in\mathbb Z}a_j2^{-j/p}
        \chi_{E_j}(x)
    \end{equation}
    holds, for some $C>0$.
\end{lemma}

It is possible to see that the best constant $C$ in \eqref{eq.lorentz.atomic}
is equivalent to the Lorentz norm of the function $f$.

The most useful properties of Lorentz spaces are the  H\"older and Young
inequalities, which extend the classical ones for Lebesgue spaces.
These were originally proved by O'Neill in \cite{oneill}.
We collect them in the following theorems:
\begin{theorem}[H\"older inequality]\label{th.holderlorentz}
    Let $f\in L^{p_1,q_1},\ g\in L^{p_2,q_2}$. The following
    estimates hold:
    \begin{itemize}
    \item if $p_1,p_2,p\in]1,\infty[$, $q_1,q_2,q\in[1,\infty]$, then
    \begin{equation}\label{eq.holderlorentz1}
        \|fg\|_{L^{p,q}}\leq
        C\|f\|_{L^{p_1,q_1}}\|g\|_{L^{p_2,q_2}},
        \qquad 1>p_1^{-1}+p_2^{-1}=p^{-1}, q_1^{-1}+q_2^{-1}
        \geq q^{-1};
    \end{equation}
    \item if $p_1,p_2\in[1,\infty[$, $q_1,q_2\in[1,\infty]$, then
    \begin{equation}\label{eq.holderlorentz2}
        \|fg\|_{L^1}\leq
        C\|f\|_{L^{p_1,q_1}}\|g\|_{L^{p_2,q_2}},
        \qquad p_1^{-1}+p_2^{-1}=1,\ q_1^{-1}+q_2^{-1}\geq1.
    \end{equation}
    \end{itemize}
\end{theorem}

We remark that the above statement does not cover the trivial inequality
\begin{equation}\label{eq.trivial}
    \|fg\|_{L^{p,q}}\leq \|f\|_{L^{\infty}}\|g\|_{L^{p,q}}
\end{equation}
which is easily proved to be true for all cases when $L^{p,q}$ is defined.

\begin{theorem}[Young inequality]\label{th.younglorentz}
    Let $f\in L^{p_1,q_1},\ g\in L^{p_2,q_2}$. Then the following
    estimates hold:
    \begin{itemize}
    \item if $p_1,p_2,p\in]1,\infty[$, $q_1,q_2,q\in[1,\infty]$,
    then
    \begin{equation}\label{eq.younglorenz1}
        \|f\ast g\|_{L^{p,q}}\leq C\|f\|_{L^{p_1,q_1}}
        \|g\|_{L^{p_2,q_2}},\qquad p_1^{-1}+p_2^{-1}
        =1+p^{-1},\ q_1^{-1}+q_2^{-1}\geq q^{-1};
    \end{equation}
    \item if $p_{1},p_2\in ]1,\infty[$,
    $q_1,q_2\in[1,\infty]$, then
    \begin{equation}\label{younglorentz2}
        \|f\ast g\|_{L^\infty}\leq C\|f\|_{L^{p_1,q_1}}
        \|g\|_{L^{p_2,q_2}},\qquad p_1^{-1}+p_2^{-1}
        =1,\ q_1^{-1}+q_2^{-1}\geq1.
    \end{equation}
    \end{itemize}
\end{theorem}

As before, we remark that the above statement does not cover
the inequality
\begin{equation}\label{eq.trivial2}
        \|f\ast g\|_{L^{p,q}}\leq C\|f\|_{L^\infty}
        \|g\|_{L^{p,q}}
\end{equation}
which is easily seen to be true in all cases when $L^{p.q}$ is
defined (e.g., by real interpolation).

We conclude this section by studying the weight
functions $w_\beta(x)=|x|(|\log|x||+1)^\beta$, with $\beta>1$
which plays a crucial role in our results; in
the following proposition we determine
precisely to which Lorentz
the powers $w_\beta^{-s}$ belong.
\begin{proposition}\label{prop.wbeta}
    For any $s>0$, $q\in[1,\infty]$ we have $w_\beta^{-s}\in L^{n/s,q}$, provided
    $\beta>1/sq$.
\end{proposition}
\begin{proof}
    We will use the equivalent Lorentz norm
    \eqref{eq.lorentz.atomic}. For any $j\in\mathbb Z$
    consider the ball $B^j:=B_{2^{j/n}}=\{x:|x|\leq 2^{j/n}\}$
    and the rings $E_j:=B^{j+1}\setminus B^j$; it is clear that
    $|E_j|=C_n2^j$, where $C_n$ depends only on
    the dimension $n$. Then, for all $x\in\mathbb R^n$ we
    have the estimate
    $$|w_\beta^{-s}(x)|=|\sum_{j\in\mathbb Z}\frac{1}
    {|x|^s(|\log|x||+1)^{\beta s}}\chi_{E_j}(x)|\leq
    \sum_{j\in\mathbb Z}(|j|\log2+1)^{-\beta s}2^{-js/n}
    \chi_{E_j}(x).$$
    The proof is concluded by the remark that the sequence
    $a_j=(|j|\log2+1)^{-\beta s}$ is in $l^q$ if and only if
    $\beta>1/sq$.
\end{proof}


\bibliographystyle{plain}

\begin{thebibliography}{10}

\bibitem{agmon}
{\sc S.~Agmon}, {\it Spectral properties of
    Schr\"odinger operators and Scattering Theory},
  Ann. Sc. Norm. Sup. Pisa Cl. Sci. {\bf 2} (1975) no. 2, 151--218.

\bibitem{ArbYaj00}
{\sc G.~Artbazar, K.~Yajima},
{\it The $L^{p}$-continuity of wave operators for one dimensional
Schr\"odinger operators},
J. Math. Sci. Univ. Tokyo {\bf 7} (2000), 221--240.

\bibitem{BL}
{\sc J.~Bergh, J.~L\"ofstr\"om},
Interpolation spaces. Springer Verlag, Berlin, 1976.

\bibitem{BRV} {\sc J.~A.~Barcelo, A.~Ruiz, L.~Vega},
    {\it Weighted Estimates for the Helmholtz Equation and
    Some Applications}, J. Funct. Anal. {\bf 150} (1997) no.2, 356--382.

\bibitem{Beals94-optiml}
    {\sc M.~Beals},
    {\it Optimal $L\sp \infty$ decay for solutions to the wave equation with
    a potential},
    Comm. Partial Differential Equations, {\bf 19} (1994) no. 7-8, 1319--1369, 1994.

\bibitem{BealsStrauss93-l}
    {\sc M.~Beals, W.~Strauss},
    {\it $L\sp p$ estimates for the wave equation with a potential},
    Comm. Partial Differential Equations, {\bf 18} (1993) no. 7-8, 1365--1397, 1993.

\bibitem{Burq04}
{\sc N.~Burq},
{\it Smoothing effect for Schr\"odinger boundary value problems},
Duke Math. J. {\bf123} no. 2 (2004), 403--427.

\bibitem{BurqPlanchonStalkerTahvildar03}
    {\sc N.~Burq, F.~Planchon, J.~Stalker, S.~Tahvildar-Zadeh},
    {\it Strichartz estimates for the wave and Schr\"odinger equations with the
    inverse-square potential},
    J. Funct. Anal. {\bf 203} (2003) no. 2, 519--549.

\bibitem{cald}
{\sc A.P.~Calder\'on},
{\it Spaces between $L^{1}$ and $L^{\infty}$ and the theorem of
Marcinkiewicz},
Studia Math. {\bf 26} (1966), 273--299.

\bibitem{Cuccagna00}
    {\sc S.~Cuccagna},
    {\it On the wave equation with a potential},
    Comm. Partial Differential Equations {\bf 25} (2000) no. 7--8, 1549-1565.

\bibitem{CuccagnaSchirmer01}
    {\sc S.~Cuccagna, P.~Schirmer},
    {\it On the wave equation with a magnetic potential},
    Comm. Pure Appl. Math. {\bf 54} (2001) no. 2, 135--152.

\bibitem{pda-vf} {\sc P.~D'Ancona, V.~Pierfelice}, {\it On the wave
    equation with a large rough potential}, submitted.

\bibitem{GeorgievVisciglia03}
    {\sc V.~Georgiev, N.~Visciglia},
    {\it Decay estimates for the wave equation with potential},
    Comm. Partial Differential Equations {\bf 28} (2003) no. 7-8, 1325--1369.

\bibitem{GinibreVelo95-generstric}
    {\sc J.~Ginibre, G.~Velo},
    {\it Generalized Strichartz inequalities for the wave equation},
    J. Funct. Anal. {\bf 133} (1995) no. 1, 50--68.

\bibitem{Goldberg04}
    {\sc M.~Goldberg},
    {\it Dispersive estimates for the three-dimensional Schr\"{o}dinger
    equation with rough potentials}, arxiv.

\bibitem{GoldbergSchlag04}
    {\sc M.~Goldberg, W.~Schlag},
    {\it Dispersive estimates for Schr\"odinger operators in dimensions one
    and three},
    Comm. Math. Phys. {\bf 251} (2004) no. 1, 157--178.

\bibitem{HassellTaoWunsch04}
    {\sc A.~Hassell, T.~Tao, J.~Wunsch},
    {\it Sharp Strichartz estimates on non-trapping asymptotically conic
    manifolds}, arxiv.

\bibitem{hormander} {\sc L.~H\"ormander}, {\rm The analysis of linear
    partial differential operators II}. Grundlehren der Mathematischen
  Wissenschaften {\bf 257}, Springer Verlag, Berlin 1983.

\bibitem{JourneSofferSogge91-decayschroed}
    {\sc J.-L.~Journ{\'e}, A.~Soffer, C.~D.~Sogge},
    {\it Decay estimates for Schr\"odinger operators},
    Comm. Pure Appl. Math. {\bf 44} (1991) no. 5, 573--604.

\bibitem{KajitaniSatoh04}
    {\sc K.~Kajitani, A.~Satoh},
    {\it Time decay estimates for linear symmetric hyperbolic systems with
    variable coefficients and its
    applications}, preprint 2004.

\bibitem{KP} 
{\sc T.~Kato, G.~Ponce},
 {\it Commutator estimates and the
Euler and Navier - Stokes equations}
Comm. Pure Appl. Math. {\bf41} (1988), 891--907.

\bibitem{KeelTao98-endpoinstric}
    {\sc M.~Keel, T.~Tao},
    {\it Endpoint Strichartz estimates},
    Amer. J. Math. {\bf 120} (1998) no. 5, 955--980.

\bibitem{LandauLifshitz}
    {\sc L.~D.~Landau, E.~M.~Lifshitz}, {\rm Course of theoretical
    physics},
    \newblock Pergamon Press, 1987.

\bibitem{Liess91-crystal}
    {\sc O.~Liess},
    {\it Decay estimates for the solutions of the system of crystal
    optics},
    Asymptotic Analysis {\bf 4} (1991) no. 1, 61--95.

\bibitem{LucenteZiliotti00-maxw}
    {\sc S.~Lucente, G.~Ziliotti},
    {\it Global existence for a quasilinear maxwell system},
    Rend. Istit. Mat. Univ. Trieste {\bf 31} (2000) suppl. 2,
    169--187.

\bibitem{MP}
{\sc J.~Mal\'y,  L.~Pick},
{\it An elementary proof of sharp Sobolev embeddings},
Proc. Amer. Math. Soc.
{\bf 130} (2002) no. 2, 555--563

\bibitem{oneill} {\sc R.~O'Neil}, {\it Convolution operators and
$L(p,q)$ spaces}, Duke Math. J. {\bf 30} (1963), 129--142.

\bibitem{PSU}
    {\sc C.~Pladdy, Y.~Saito, T.~Umeda}, {\it Resolvent estimates of
    the Dirac operator}, Analysis {\bf 15} (1995) no. 2, 123--149.

\bibitem{reed-simonI}{\sc M.~Reed, B.~Simon}, {\rm Methods of Modern
    Mathematical Physics vol. I: Functional Analysis}. Academic Press,
    New York, San Francisco, London
    1980.

\bibitem{reed-simonII} {\sc M.~Reed, B.~Simon}, {\rm Methods of Modern
    Mathematical Physics vol. II: Fourier Analysis,
    Self-Adjointness}. Academic Press, New York, San Francisco, London
  1975.

\bibitem{RobbianoZuily05}
    {\sc L.~Robbiano, C.~Zuily},
    {\it Strichartz estimates for Schr\"odinger equations with variable
    coefficients}, arxiv.

\bibitem{RodnianskiSchlag04-disp}
    {\sc I.~Rodnianski, W.~Schlag},
    {\it Time decay for solutions of Schr\"odinger equations with rough and
    time-dependent potentials},
    Invent. Math. {\bf 155} (2004) no. 3, 451--513.

\bibitem{Schlag05}
    {\sc W.~Schlag},
    {\it Dispersive estimates for Schr\"odinger operators: A
    survey}, preprint 2005.

\bibitem{SS}
    {\sc J.~Shatah, M.~Struwe},
    {\it Geometric wave equations}.
    New York University Courant Institute of Mathematical Sciences,
    New York, 1998.

\bibitem{StaffilaniTataru02-stricschroed}
    {\sc G.~Staffilani, D.~Tataru},
    {\it Strichartz estimates for a Schr\"odinger operator with nonsmooth
    coefficients},
    Comm. Partial Differential Equations {\bf 27} (2002) no. 7-8,
    1337--1372.

\bibitem{Stein} {\sc E.~Stein, R.~Shakarchi}, {\rm Complex
    analysis}, Princeton University Press, Princeton 2003.

\bibitem{Tarulli04-smoo}
    {\sc M. Tarulli Di Giallonardo},
    {\it Resolvent estimates for scalar field with electromagnetic
    perturbation}, El. J. Diff. Eqns {\bf 2004} (2004), 1-18.

\bibitem{taylor}
    {\sc M. Taylor},
    {\it Partial differential equations, Vol. II}.
     Springer, New York 1996.

\bibitem{Yajima95-waveopN}
    {\sc K. Yajima},
    {\it The $W^{k,p}$-continuity of wave operators for Schr\"odinger
    operators},
    J. Math. Soc. Japan {\bf 47} (1995) no. 3, 551--581.

\bibitem{Yajima95-waveopNeven}
    {\sc K. Yajima},
    {\it The $W^{k,p}$-continuity of wave operators for Schr\"odinger
    operators III, even dimensional cases $m\geq4$},
    J. Math. Sci. Univ. Tokyo {\bf 2} (1995), 311--346.

\bibitem{Yajima99-waveop2D}
    {\sc K. Yajima},
    {\it $L\sp p$-boundedness of wave operators for two-dimensional
    Schr\"odinger operators},
    Comm. Math. Phys. {\bf 208} (1999) no. 1, 125--152.

\bibitem{yamada} {\sc O.~Yamada}, {\it On the principle of limiting
    absorption for the Dirac operators},
  Publ. Res. Inst. Math. Sci. {\bf 8} (1972/73), 557--577.

\end{thebibliography}

\end{document}